\documentclass[10pt,a4paper]{article}

\usepackage[authoryear]{natbib}

\usepackage[english]{babel}
\usepackage[utf8]{inputenc}
\usepackage{amsmath,amssymb,amsthm}
\usepackage[dvipsnames]{xcolor}
\usepackage{float}
\usepackage{bm}
\usepackage{booktabs,enumerate,fullpage}
\usepackage{rotating}
\usepackage{url}
\usepackage{multirow}
\usepackage{mathtools}
\usepackage{graphicx}
\usepackage{MnSymbol}
\usepackage[margin=3cm]{geometry}

\usepackage[ruled,linesnumbered]{algorithm2e}
\usepackage{comment}
\allowdisplaybreaks

\newtheorem{theorem}{Theorem}
\newtheorem{lemma}{Lemma}

\begin{document}

\title{Chemotherapy planning and multi-appointment scheduling: formulations, heuristics and bounds}

\author{
	Giuliana Carello\thanks{Politecnico di Milano, Department of Electronics, Information and Bioengineering, piazza Leonardo da Vinci 32, 20133 Milan, Italy}
	\and 
	Mauro Passacantando\thanks{University of Milano-Bicocca, Department of Business and Law, via Bicocca degli Arcimboldi 8, 20126 Milan, Italy} 
	\and 
	Elena Tanfani\thanks{University of Genova, Department of Economics, via Vivaldi 5, 16126 Genova, Italy}
	}

\date{}

\maketitle

\noindent\textbf{Abstract.}
The number of new cancer cases is expected to increase by about 50\% in the next 20 years, and the need for chemotherapy treatments will increase accordingly. 
Chemotherapy treatments are usually performed in outpatient cancer centers where patients affected by different types of tumors are treated. The treatment delivery must be carefully planned to optimize the use of limited resources, such as drugs, medical and nursing staff, consultation and exam rooms, and chairs and beds for the drug infusion. 
Planning and scheduling chemotherapy treatments involve different problems at different decision levels. 
In this work, we focus on the patient chemotherapy multi-appointment planning and scheduling problem at an operational level, namely the problem of determining the day and starting time of the oncologist visit and drug infusion for a set of patients to be scheduled along a short-term planning horizon. 
We use a per-pathology paradigm, where the days of the week in which patients can be treated, depending on their pathology, are known. 
We consider different metrics and formulate the problem as a multi-objective optimization problem tackled by sequentially solving three problems in a lexicographic multi-objective fashion. 
The ultimate aim is to minimize the patient's discomfort. The problems turn out to be computationally challenging, thus we propose bounds and ad-hoc approaches, exploiting alternative problem formulations, decomposition, and $k$-opt search. The approaches are tested on real data from an Italian outpatient cancer center and outperform state-of-the-art solvers. 

\


\section{Introduction} 
\label{sec:intro}

In the last decades, we observed an increase in cancer incidence in most of the Western and industrialized countries. According to a study by the International Agency for Research on Cancer on the impact of population aging on the future cancer burden, the global cancer burden is expected to be 28.4 million cases in 2040, which corresponds to a 47\% rise since 2020~\citep{SFSLS2021}. 
The World Health Organisation identified cancer as the second leading cause of death globally, with 10 million deaths worldwide in 2020 and accounting for one in six deaths worldwide~\citep{WHO2020}. 
However, the overall cancer mortality rate has been steadily decreasing over the last 30 years~\citep{SMWJ2023}. 
This decrease can be related to many factors, e.g., the advancement of prevention, screening, and early detection strategies, the addition of chemotherapy treatment after surgery for some cancer types, combined therapies for many cancers, and survivorship care~\citep{SFSLS2021}.

Features such as the type of cancer, its site, the stage of illness and severity, guide to select the treatment options and their combination. 
The most widely used traditional treatment methods are surgery, chemotherapy, and radiotherapy. Recently, modern modalities that include hormone therapy, anti-angiogenic, stem cell therapies, immunotherapy, and dendritic cell-based immunotherapy are receiving growing attention~\citep{DMH2021}. 
However, currently, chemotherapy is still the most used evidence-based treatment for the majority of cancer types due to its proven ability to reduce morbidity and mortality~\citep{ANAND2023}. 

Recent studies estimated that between 2018 and 2040, the number of patients requiring first-course chemotherapy will annually increase from 9.8 million to 15.0 million, with a relative increase of 53\%, and predicted a major health crisis driven by unmet chemotherapy demand~\citep{WLS2019}. 
This dramatically increased demand forces private and public healthcare organizations to design improved solutions. 
Congestion and workload issues have been reported in many studies, and improved scheduling practices have been identified, as a key component to mitigate resource shortage and ensure timely and effective cancer treatment~\citep{ILK2021,HEDEDK2020}.  Indeed, poor organization of the outpatient chemotherapy process creates long waiting times for the infusion and patient dissatisfaction~\citep{MHT2020}. 

Many healthcare organizations are reorganizing the delivery of chemotherapy treatments into dedicated facilities, referred to as outpatient cancer centers or infusion centers, to improve resource use and patient care~\citep{VHARTEN2014}. 
This reorganization process changed the delivery paradigm from a hospital specialty-based to a patient pathology-based approach and poses problems and challenges. 

Chemotherapy appointment planning and scheduling in outpatient settings involves different levels of decisions hierarchically linked~\citep{LAJOST2016}. 
At a strategic level, the main decision pertains to the dimensioning of the outpatient center with respect to the human and physical resources involved in the chemotherapy delivery process (e.g., nurses, doctors, pharmacists, paramedical and technical staff, beds, chairs, exam and consultation rooms, etc.). 
At a tactical planning level, the available resources must be assigned to a set of macro groups of cancer pathologies. In particular, the Master Chemotherapy Planning (MCP) gives the days when the different cancer pathologies are scheduled over a mid-term planning horizon, together with the clinician rostering that covers the cyclic schedule based on clinicians availability and skills~\citep{CLTT2022}. 
Then, at an operational level, the scheduling of patients' appointments must be determined along a shorter planning horizon (usually a week). 
The problem is divided into two subproblems: determining the day of the infusion (interday scheduling) and the time slot of the appointment (intraday scheduling)~\citep{CSL2023}.

In this work, we focus on the operational level of a chemotherapy appointment scheduling problem in a multi-stage and multi-resource environment. 
Following the same-day outpatient chemotherapy policy~\citep{HEP2022}, patients undergo a blood test, then they are visited by an oncologist who, depending on the overall clinical status and health condition of the patients, decides if they can receive the chemotherapy treatment in the infusion room or not. 
The process takes place in two consecutive days in case of next-day policy.

We tackle the problem of determining the appointment date and starting time of the oncologist visit and the drug infusion for a set of patients over a one-week planning horizon jointly addressing the interday and intraday scheduling problem. 
The problem belongs to the class of multi-appointment scheduling problems~\citep{MD2020}. 

We formulate the problem as a multi-objective optimization problem aimed at optimizing different criteria and metrics of interest in a patient-centered perspective, i.e., number of scheduled patients, waiting time of patients in the center, and patient preferences on the infusion place (beds or chairs). 

Maximizing the number of treated patients can be seen as a feasibility check of the allocated resources. 
In fact, treatment plans are made by oncologists for each patient according to existing chemotherapy protocols based on clinical trials. 
As a result, chemotherapy treatments are given in cycles with rest periods during and between each cycle. The number and duration of the cycles change according to the cancer type and stage of the disease, and it is important to strictly adhere to the patient treatment plan to achieve the best health outcomes. 
Waiting time is an important quality indicator that affects the level of satisfaction of the patients, enhancing the patient treatment experience, and ultimately improving patient care and outcomes~\citep{RK2017}. 
Finally, focusing on a patient-centered perspective the last metric considers the patient preferences and satisfaction about the place of the drug administration (chair or bed). 
Obviously, critical patients must receive the infusion in a bed due to their clinical status, whereas non-critical patients prefer to have the drug administered in a chair, where they can read a book or chat with other patients. 

As the objectives are not equally important, we model the problem using a lexicographic approach and we solve it by optimizing a sequence of Integer Linear Programming (ILP) models. 
 The problems turn out to be computationally challenging, thus
alternative formulations are proposed and proved to be equivalent. 
Ad-hoc solution approaches for each problem and bounds are proposed to speed up the solution process.

The organization of the paper is the following. 
In Section~\ref{sec:litRew}, the relevant literature on patient appointment planning and scheduling for chemotherapy treatments is analyzed. 
The problem and the multi-objective formulation are described in Section~\ref{sec:problem}. 
The solution approaches and the bounds proposed
for each problem are reported and explained in Section~\ref{sec:procedure}. 
Results are discussed in Section~\ref{sec:results} and conclusions and further developments are given in Section~\ref{sec:conclusions}.


\section{Literature review and contribution of the paper} 
\label{sec:litRew}

Chemotherapy planning and scheduling encompass many challenges and problems, such as determining optimal drug combinations, dosages, and treatment cycles, planning treatments, scheduling patient appointments, and allocating the main resources involved with the care process~\citep{SSB2019}.  
In the past, the attention of researchers has been mainly focused on the medical issues of the problem. 
Many papers have been published proposing methods to deal with the optimal design of cancer chemotherapy plans balancing the benefits of treating tumors with the adverse toxic side effects caused by the anti-cancer drugs, taking into account several factors, such as tumor growth features and drug infusion rates~\citep{AHL2006,SAES2014,RH2015}.

More recently, the scientific literature has been also directed toward operations management problems.
In~\cite{LAJOST2016}, the first systematic literature review on outpatient chemotherapy planning is reported. The review classifies the papers based on several features, such as the planning level (strategic, tactical, and operational), the objective of the study and related performance metrics, the type of resources, and the uncertainty sources considered. 
An updated optimization-oriented systematic literature review in the period 2009-2021 is proposed in~\cite{HEM2022}. 
The authors first classified the contributions into three classes depending on the main problem addressed, i.e. planning, scheduling, and assignment problems. 
The 45 publications selected have been then grouped with respect to the scope of the analysis, the considered complexities, the modeling techniques, and the solution methods proposed. 
Planning problems include determining treatment days, recovery days, and the drug doses used each day (treatment scope), determining the days of the infusion according to the treatment plans (patient scope), and determining the working periods of the oncologists and balancing the workload of the medical staff (oncologist scope). 
Scheduling problems are mainly focused on patient scheduling on the day of the infusion, nurse rostering and drug preparation and  delivery~\citep{KGB2017}. 

In line with the above referenced classifications, we focus the analysis of the literature on the studies published in international journals in the period 2019-2023 dealing with chemotherapy appointment planning and scheduling problems, focusing on a patient scope and operational level. Note that, in some publications, such planning and scheduling problems are also referred to as interday and intraday appointment scheduling~\citep{CSL2023}.

The problem of determining the appointment date of arriving patients within a specific time window for treatment has been tackled in~\cite{GP2014}. 
The authors propose a Markov Decision process formulation and apply linear-programming-based Approximate Dynamic Programming to obtain approximate solutions. 
The same problem has been faced in~\cite{HEDEDK2020}, using a multi-criteria Mixed-Integer Linear Programming (MILP) model, based on a three-stage heuristic. 
The aim is to minimize the number of deferring times for scheduled appointments taking into account nurse capacity constraints. 

Several contributions tackled the patient scheduling problem proposing exact formulations, heuristics, and meta-heuristics~\citep{HEM2022}. 
Many papers jointly considered patients-to-resources assignment decisions, such as nurses, oncologists, beds and chairs, and pharmacists. 
In~\cite{HDK2019} the authors developed an ILP model in order to assign the appointment slots to patients subject to the nursing constraints. 
They use a parametric multi-objective function for minimizing the weighted flowtime and minimizing the makespan. 
~\cite{AHDK2020} proposed a multi-criteria MILP model to schedule patient appointments with starting times close to patients’ ready times, balanced workload among nurses, few nurse changes during appointments, and few nurse full-time equivalents assigned to the schedule of the day in order to minimize the labor cost.  
In~\cite{LCA2023} an ILP model is proposed to schedule all the patients assigned to a short planning horizon. 
They also account for the preparation time of the medication to be administered. 
They proposed a linear relaxation of the model, based on treatment patterns, that is solved via column generation. 
The same-day chemotherapy outpatient scheduling problem has been tackled in~\cite{CCF2022} as a multi-stage problem with the aim of reducing idle times within and among the different stages, namely medical consultation, drug preparation, and infusion. \cite{ABI2019} model the multi-stage appointment problem as a deterministic flexible job shop in which patients correspond to jobs and activities to the stages of the flow of the patients within the cancer clinic. 
The objective of the model is the minimization of the total amount of time spent by patients on the day of the treatment.

Different sources of uncertainty have been considered in the literature.
In~\cite{ILK2021} a MILP robust slack allocation model is designed to find the optimal patient appointment schedule, considering the fact that the infusion time of patients may take longer than expected. 
This model minimizes the weighted sum of the total waiting times of patients, the makespan, and the number of beds used through the planning horizon. 
A robust scheduling heuristic is developed based on the adaptive large neighborhood search to determine patient appointments of real-size infusion centers. 
The problem of sequencing the patients for a one-day session of chemotherapy infusion considers the uncertainty related to the patient's inability to receive the treatment is solved in~\cite{GAROXI2020} using a stochastic optimization model. 
A slightly different version of the problem is considered in~\cite{KGC2022}, where the problem of scheduling patient appointments and assigning patients to nurses considering the uncertainty in infusion durations in a given day is addressed. 
The authors formulate a two-stage stochastic MILP model with the goal of minimizing the expected weighted sum of excess patient acuity, waiting time, and nurse overtime. 
An online optimization approach is proposed in~\cite{HTC2020} to deal with patients’ multiple requests for chemotherapy treatments. 
The authors propose an adaptive and flexible procedure that combines two optimization models. The first one schedules incoming appointment requests, while the second one reschedules already booked appointments with the goal of better allocating resources as new information becomes available. 

Few integrated approaches that combine planning and scheduling problems have been recently proposed. 
In~\cite{HEEL2019} a two-phase approach is proposed. First, a MILP model is proposed for determining the day of treatment for a set of patients and finding the optimal number of needed nurses and pharmacists. 
Then, in the second phase, a discrete event simulation model is used to generate patient appointment schedules that minimize the treatment delay for patients and the total completion times of treatments each day.  Similarly, in~\cite{BELARO2019} a two-stage procedure for chemotherapy appointment scheduling of new patients and for assigning the daily patient mix to available nurses is introduced. 
In the first phase, the appointments of the new patients are assigned at the end of each day, and the daily requirement for nurses is determined. 
In the second phase, patients are assigned to nurses while minimizing the number of nurses required. 
In~\cite{CSL2023} a two-phase solution approach is developed to deal with the interday and intraday scheduling determining the date and time slot of the infusion of a set of patients to be scheduled. 
The two subproblems have been solved sequentially with optimization models developed for each phase. The solutions have been evaluated in terms of makespan, resource utilization, overtime, and patient diversion metrics.
~\cite{HDK2023} use a fixed template of slots for the online scheduling of appointments simultaneously addressing the interday and intraday scheduling problem.

The analysis of the literature shows that, although some integrated approaches for interday and intraday scheduling have been proposed, none of them tackled the problem as a multi-appointment problem. 
Only one paper~\citep{CCF2022} deal with the multi-appointment scheduling problem considering more than one activity to be executed in sequence in a given planning horizon. However, they do not tackle simultaneously the problem of determining the date of the appointment. 

As for the objective functions considered in the studies, time and cost metrics received
more attention than workload and satisfaction measures~\citep{HEM2022}.
Only one paper seeks to maximize patient preferences and to minimize the number of unscheduled patients with respect to the requests~\citep{HTC2020}.  

The novelty of this paper is twofold. 
First, we jointly address the interday and intraday multi-appointment scheduling and bed/chair assignment decisions. 
As far as we know, this is the first attempt to integrate patient-centered operational decisions in a unified approach that follows the trajectory of the patients in the cancer center during the day of the treatment. 
The integrated planning and scheduling solution determines simultaneously, for a short-term planning horizon (usually a week), the day of the appointment (planning decision) as well as the starting time of the oncologist visit and of the infusion (scheduling decision) for a set of patients to be scheduled according to the treatment plans. 
The assignment of patients to the available beds and chairs is also determined. 
Secondly, following the patient-centered scope, we combine different metrics already proposed in previous studies, such as minimizing the number of unscheduled appointment requests, minimizing patient waiting time, and maximizing patient preferences, into a multi-objective optimization approach, where the ultimate aim is to minimize the patient discomfort. 
The clinical consistency of the treatment plans is ensured as the set of patients to be treated each week has been determined based on clinical analysis, and the first and most relevant metric aims to ensure that they are all scheduled. 
The clinicians covering of the schedule is ensured by taking as data input the solution of a MILP model proposed in a recently published work~\citep{CLTT2022}. 


\section{Problem description and formulation} \label{sec:problem}

We tackle the problem of determining the appointment date and starting time of the consultation visit and chemotherapy infusion for a set of patients $P$.  
Patients are classified into a set $K$ of cancer pathologies: binary parameter $u_{pk}$ is equal to 1 if the main pathology of patient $p$ is $k$, and equal to 0 otherwise. 
We consider a planning horizon represented by a set $T$, usually one working week long. 
Each working day $t \in T$ is divided into a set $H$ of time slots.
Each patient must undergo first the clinician consultation (visit) and then the chemotherapy infusion, on the same day. 
For each patient $p$, the consultation visit time $v_p$ and the infusion time $f_p$ are given, and assumed to be deterministic values. 
Consultations are carried on in exam rooms. 
The set of exam rooms is denoted as $R$, and the subset of time slots when the exam rooms are open on each working day is denoted as $H_V$. 
We assume that each exam room is devoted to a single pathology on each day: let the binary parameter $w_{rkt}$ be equal to 1 if room $r$ is assigned to pathology $k$ in day $t$, and equal to 0 otherwise, according to the given MCP obtained as a result of the tactical level decision~\citep{CLTT2022}.
As for the infusion, we can use both chairs and beds depending on the clinical status of the patients. 
The set $P$ of patients is partitioned into two subsets: $P_B$ includes the patients who must receive the infusion in a bed (\textit{critical patients}), while the patients in $P \setminus P_B$ can receive the infusion both in a chair or in a bed (\textit{non-critical patients}). 
The sets of beds and chairs are denoted with $B$ and $S$, respectively. 
Note that the chair is always the most comfortable and preferable choice for non-critical patients.

We assume that the laboratory for the blood tests as well as the pharmacy for the drug preparation are within the outpatient center and always available to meet the requests. 
Similarly, we assume that the nurses are always available to follow the clinical trajectory of all patients scheduled each day (see~\cite{GAROXI2020} for similar assumptions). 

In our approach, the resources that limit the patient flow are the exam rooms (and the assigned oncologists), the beds, and the chairs. 
Due to their capacity, it may not be always possible to treat all the patients in the considered time horizon according to the treatment plan. 
Therefore, we have to decide which patients are treated in the considered time horizon, and the starting time of consultation visit and chemotherapy infusion for the selected patients. 
However, maximizing the number of scheduled patients is not the only goal of the problem. We also want to minimize the waiting time between the visit and the infusion, thus minimizing the time spent by patients in the center. 
In order to minimize patients' discomfort we also want to assign as many non-critical patients as possible to chairs. 
The three goals are clearly not equally relevant, so the problem is tackled as a multi-objective problem with a lexicographic objective function. 

The problem belongs to the family of multi-appointment scheduling problems~\citep{MD2020}. 
It shares some features with the flexible flow shop, with patients acting as jobs. 
However, it is more complex than the flexible flow shop. Indeed, differently from the classical flexible flow shop problem, patients are partitioned into subsets. 
For the first stage (visit), patients are divided based on their pathology, and for the second stage (infusion) they are divided into critical or non-critical patients. 
A patient in a set cannot be processed by any machine, but only by a subset of machines, both for the visit and for the infusion, and patients of different subsets of the first stage shared resources and machines in the second one. 

The problem can be proved NP-hard. 
It is possible to reduce the problem where only the first objective -- the maximization of the number of treated patients -- is considered and only the infusion on chairs from the bin packing problem, where the patients represent the items, and the chairs represent the bins.

The problem can be formulated as described in the next sections.


\subsection{Formulation: variables and constraints}
\label{subsec:variables-and-const}

The problem can be formulated with four families of binary variables.
The first one represents the day $t$, starting time $h$, and exam room $r$ of the consultation visit of patient $p$:
$$
\alpha_{pthr} = 
\begin{cases}  
1 & \text{if patient $p \in P$ starts the visit in room $r$, in day $t$ at time slot $h$,} 
\\
0 & \mbox{otherwise.}  
\end{cases}
$$
The variables are defined so as to prevent the assignment of a patient to a day and room where a suitable clinician is not present according to the MCP, described by parameter $w_{rkt}$, and to guarantee that the visits are completed within the last available slot each day.

The second variable assigns the starting time and the bed for the infusion to a critical patient $p \in P_B$ requiring a bed:
$$
    \beta_{pthb} = 
    \begin{cases}  
    1 & \text{if patient $p \in P_B$ starts the infusion in bed $b$, in day $t$ at time slot $h$,} 
    \\
    0 & \mbox{otherwise.}  
    \end{cases}
$$
For non-critical patients, two additional variables are introduced, which determine whether the patient receives the infusion in a bed or in a chair and the starting time for the infusion:
$$
\gamma^B_{pthb}=
\begin{cases} 
1 & \text{if patient $p \in P \setminus P_B$ starts the infusion in bed $b$, in day $t$ at time slot $h$,} 
\\
0 & \mbox{otherwise,}
\end{cases}
$$
and 
$$
\gamma^S_{pths}=
\begin{cases} 
1 & \text{if patient $p \in P \setminus P_B$ starts the infusion in chair $s$, in day $t$ at time slot $h$,}
\\
0 & \mbox{otherwise.}
\end{cases}
$$

Variables $\beta_{pthb}, \gamma^B_{pthb}$ and $\gamma^S_{pths}$ are defined so as to guarantee that the infusions end before the last available time slot each day.

The constraints of the problem can be formulated as follows:
\begin{align}
\label{C:visit_assignment}
& \sum_{t \in T} \sum_{h \in H_V} \sum_{r \in R} \alpha_{pthr} \leq 1 
& \forall p \in P 
\\
\label{C:onePerRoom}
& \sum_{p \in P} \sum_{q = \max\{1,h+1-v_p\}}^h \alpha_{ptqr} \leq 1 
& \forall t \in T, h \in H_V, r \in R 
\\
\label{C:visit_infusion_1}
& \sum_{h \in H_V}\sum_{r \in R} \alpha_{pthr} 
= \sum_{h \in H} \sum_{b \in B} \beta_{pthb} 
& \forall\ t \in T, p \in P_B
\\
\label{C:visit_infusion_0}
& \sum_{h \in H_V}\sum_{r \in R} \alpha_{pthr} 
= \sum_{h \in H} \sum_{b \in B} \gamma^B_{pthb} 
    + \sum_{h \in H} \sum_{s \in S} \gamma^S_{pths} 
& \forall t \in T, p \in P \setminus P_B
\\
\label{C:first_visit_then_infusion_1}
& \sum_{r \in R} \sum_{h \in H_V} (h+v_p) \alpha_{pthr} 
    \leq \sum_{h \in H} \sum_{b \in B} h \, \beta_{pthb} 
& \forall t \in T, p \in P_B
\\
\label{C:first_visit_then_infusion_0}
& \sum_{r \in R} \sum_{h \in H_V} (h+v_p) \alpha_{pthr} 
    \leq \sum_{h \in H} h \, \left[ \sum_{b \in B} \gamma^B_{pthb}
         + \sum_{s \in S} \gamma^S_{pths} \right]
& \forall t \in T, p \in P \setminus P_B
\\
\label{C:capacitychairs}
& \sum_{p \in P \setminus P_B} \sum_{q = \max\{1, h+1-f_p\}}^h \gamma^S_{ptqs} 
\leq 1 
& \forall t \in T, h \in H, s \in S
\\
\label{C:capacityBeds}
& \sum_{p \in P_B} \sum_{q = \max\{1, h+1-f_p\}}^h \beta_{ptqb}
   + \sum_{p \in P \setminus P_B} \sum_{q = \max\{1, h+1-f_p\}}^h \gamma^B_{ptqb}
\leq 1 
& \forall t \in T, h \in H, b \in B.
\end{align}
Constraints~\eqref{C:visit_assignment} guarantee that each patient is visited at most once in the planning horizon. 
Constraints~\eqref{C:onePerRoom} guarantee that at most one patient is visited in each time slot in each exam room. 
Constraints~\eqref{C:visit_infusion_1} and \eqref{C:visit_infusion_0} guarantee that the consultation visit and the infusion take place on the same day, for critical ($P_B$) and not-critical ($P \setminus P_B$) patients, respectively, while constraints \eqref{C:first_visit_then_infusion_1} and \eqref{C:first_visit_then_infusion_0} guarantee that the infusion starts only after the consultation visit is finished. Constraints~\eqref{C:capacitychairs} prevent more than one patient from receiving the infusion on the same day and time slot in the same chair.
Similarly, constraints~\eqref{C:capacityBeds} limit to one the number of patients who  receive the infusion in a bed, for each day and time slot.


\subsection{Formulation: goals and metrics}
\label{hierarchicalObj}

The first goal is maximizing the number of patients treated in the planning horizon. 
The metrics can be formulated as follows:
\begin{eqnarray}
\label{OBJ:maxVisit}
\varphi_1(\alpha) = \sum_{p \in P} \sum_{t \in T} \sum_{h \in H_V} \sum_{r \in R} \alpha_{pthr}.
\end{eqnarray}
The second goal is to minimize the time spent by patients in the center. 
Since the visit and infusion time are given for each patient $p$, the second goal results in minimizing the idle time of patients, namely the difference between the starting time of the infusion and the end of the visit.
In order to obtain balanced metrics for this goal, we decided to focus on the maximum waiting time among all patients each day. 
We define $\varphi_A^t$ as the maximum waiting time in day $t$ for the critical patients: 
\begin{eqnarray}
\varphi_A^t=    \max_{p \in P_B} \left\{ 
 \sum_{h \in H} \sum_{b \in B} h \beta_{pthb} 
- 
\sum_{h \in H_V} \sum_{r \in R} (h+v_p) \alpha_{pthr}
\right\}.
\end{eqnarray}
Analogously, we define $\varphi_B^t$ as the maximum waiting time in day $t$ for the non-critical patients:
\begin{eqnarray} 
\varphi_B^t= \max_{p \in P \setminus P_B} \left\{ 
\sum_{h \in H} h \left[
\sum_{b \in B} \gamma^B_{pthb} + \sum_{s \in S} \gamma^S_{pths}
\right]
- 
\sum_{h \in H_V} \sum_{r \in R} (h+v_p) \alpha_{pthr}
\right\}.
\end{eqnarray}
The maximum waiting time $\varphi^t_2$ in day $t$ among all patients is thus computed as follows:
\begin{equation*}
\varphi^t_2(\alpha,\beta,\gamma^B,\gamma^S)=
\max \left\{ \varphi_A^t, \varphi_B^t \right\}.
\end{equation*}
Then, we sum such maximum waiting times over the working days of the planning horizon. The resulting metrics to be minimized is formulated as follows: 
\begin{equation}
\label{OBJ:waitingTime}
\varphi_2(\alpha,\beta,\gamma^B,\gamma^S) = 
\sum_{t \in T}
\varphi^t_2(\alpha,\beta,\gamma^B,\gamma^S).
\end{equation}
Finally, we take into account the preferences of non-critical patients: although they can receive the infusion in a bed or in a chair, the chair is the most comfortable and preferable choice. 
The corresponding metrics to be maximized can be represented as follows:
\begin{eqnarray} 
\label{OBJ:chairPreferences}
    \varphi_3(\gamma^S) = 
    \sum_{p \in P \setminus P_B} \sum_{t \in T} \sum_{h \in H} \sum_{s \in S} \gamma^S_{pths}.
\end{eqnarray}

As mentioned, the three metrics are not equally relevant: maximizing $\varphi_1$
is more important than minimizing $\varphi_2$,
which, in turn, is more important than maximizing $\varphi_3$. 


\section{Solution approach} 
\label{sec:procedure}

We tackle the lexicographic multi-objective problem by applying an $\epsilon$-constrained approach, solving a sequence of three models. When we optimize a goal, we force the value of the previously optimized objectives to the best one. 

The first problem ($P_1$) aims to maximize the number of scheduled patients and can be modeled as follows:
\begin{align}
\label{P1}
\tag{$F_1$} 
\left\{
\begin{array}{ll}
\max\limits_{(\alpha,\beta,\gamma^B,\gamma^S)}\ & \varphi_1(\alpha) 
\\
\text{subject to} & \eqref{C:visit_assignment}-\eqref{C:capacityBeds}.
\end{array}
\right. 
\end{align}
The optimal value of $(P_1)$ is stored in a parameter $\bar{v}_1$. 

Then, we consider problem $(P_2)$ that aims to find, among all the equivalent optimal solutions of $(P_1)$, one that minimizes the sum of the daily maximum waiting times. 
Problem $(P_2)$ can be modeled as follows: 
\begin{align}
\label{P2}
\tag{$F_2$} 
\left\{
\begin{array}{ll}
\min\limits_{(\alpha,\beta,\gamma^B,\gamma^S)}\ & \varphi_2(\alpha,\beta,\gamma^B,\gamma^S) 
\\
\text{subject to} & \eqref{C:visit_assignment}-\eqref{C:capacityBeds}
\\
& \varphi_1(\alpha) \geq \bar{v}_1.
\end{array}
\right.
\end{align}
The maximum waiting time in each day $t$ corresponding to an optimal solution of $(P_2)$ is stored in a parameter $\bar{v}^t_2$. 

Finally, the last problem $(P_3)$ aims to maximize the number of non-critical patients who  receive the treatment on a chair, while keeping the number of treated patients not worse than $\bar{v}_1$ and the daily maximum waiting time not worse than $\bar{v}^t_2$ for any $t \in T$. 
Problem $(P_3)$ can be modeled as follows:
\begin{align}
\label{P3}
\tag{$F_3$} 
\left\{
\begin{array}{rll}
\max\limits_{(\alpha,\beta,\gamma^B,\gamma^S)}\ & \varphi_3(\gamma^S) & 
\\
\text{subject to} & \eqref{C:visit_assignment}-\eqref{C:capacityBeds} &
\\[2mm]
& \varphi_1(\alpha) \geq \bar{v}_1 &
\\[2mm]
& \varphi^t_2(\alpha,\beta,\gamma^B,\gamma^S) \leq \bar{v}^t_2 & \qquad \forall\ t \in T.
\end{array}
\right.
\end{align}
To speed up the solution process, problems $(P_1)$, $(P_2)$ and $(P_3)$ are not solved through the formulations \eqref{P1}, \eqref{P2} and \eqref{P3}, but we devised specific solution approaches for each of them,  described in the following subsections.


\subsection{Solving problem $(P_1)$}
\label{s:solvingP1}

Solving the formulation \eqref{P1} is very time-consuming (see Section~\ref{sec::results-P1}). 
Thus, we devised an alternative formulation, based on a smaller set of variables, and we proved that it is equivalent to~\eqref{P1}.
The alternative model is based on an aggregate formulation where, instead of considering each room separately, rooms are aggregated into a single capacitated resource, and similarly for chairs and beds. 
The aggregate formulation uses four families of binary variables. 
The first one represents the choice of the day and starting visit time of patient $p$, which can be any time slot $h$ before the last time slot available for visits $|H_V|$:
$$
    x_{pth} = 
    \begin{cases} 
    1 &\text{if patient $p \in P$ starts the visit in day $t \in T$, at time slot $h\in H_V$,} 
    \\
    0 & \text{otherwise.}  
    \end{cases}
$$ 
The second variable assigns a starting time for the infusion of a critical patient, before the last available slot for infusions $|H|$:
$$
    y_{pth} = \begin{cases}  
    1 & \text{if patient $p \in P_B$ starts the infusion (in some bed) in day $t \in T$, at time slot $h \in H$,} 
    \\
    0 & \text{otherwise.}  
    \end{cases}
$$ 
For non-critical patients, two additional variables are introduced, which tell whether the patient receives the infusion in a bed or in a chair and assign the patient's infusion starting time:
$$
    z^B_{pth}=
    \begin{cases} 
    1 & \text{if patient $p \in P \setminus P_B$ starts the infusion in a bed, in day $t \in T$, at time slot $h \in H$,} 
    \\
    0 & \text{otherwise,}
    \end{cases}
$$
and 
$$
    z^S_{pth}=
    \begin{cases} 
    1 & \text{if patient $p \in P \setminus P_B$ starts the infusion in a chair, in day $t \in T$, at time slot $h \in H$,}
    \\
    0 & \text{otherwise.}
    \end{cases}
$$
Variables $x_{pth}$, $y_{pth}$, $z^B_{pth}$ and $z^S_{pth}$ are defined so as to avoid assignments or starting times unacceptable either because the pathology is not treated in the day or because the visit/infusion would end when the center is already close.

The aggregate formulation $(AF_1)$ of problem $(P_1)$ is defined as follows:
\begin{align}
\label{obf:max_p_surr}
& \max\ F_1(x) := \sum_{p \in P} \sum_{t \in T} \sum_{h \in H_V} x_{pth} 
\end{align}
subject to:
\begin{align}
\label{Cs:visit_assignment}
& \sum_{t \in T} \sum_{h \in H_V} x_{pth} \leq 1 
& \forall p \in P 
\\
\label{Cs:onePerRoom}
& \sum_{\substack{p \in P: \\ u_{pk}=1}} \sum_{q = \max\{1,h+1-v_p\}}^h x_{ptq} 
\leq \sum_{r \in R} w_{rkt} 
& \forall t \in T, h \in H_V, k \in K
\\
\label{Cs:visit_infusion_1}
& \sum_{h \in H_V} x_{pth} 
= \sum_{h \in H} y_{pth} 
& \forall t \in T, p \in P_B
\\
\label{Cs:visit_infusion_0}
& \sum_{h \in H_V} x_{pth} 
= \sum_{h \in H} (z^B_{pth} + z^S_{pth}) 
& \forall t \in T, p \in P \setminus P_B
\\
\label{Cs:first_visit_then_infusion_1}
& \sum_{h \in H_V} (h+v_p) x_{pth} 
    \leq \sum_{h \in H} h y_{pth} 
& \forall t \in T, p \in P_B
\\
\label{Cs:first_visit_then_infusion_0}
& \sum_{h \in H_V} (h+v_p) x_{pth} 
    \leq \sum_{h \in H} h (z^B_{pth} + z^S_{pth})
& \forall t \in T, p \in P \setminus P_B
\\
\label{Cs:capacitychairs}
& \sum_{p \in P \setminus P_B} \sum_{q = \max\{1, h+1-f_p\}}^h z^S_{ptq} 
\leq |S| 
& \forall t \in T, h \in H
\\
\label{Cs:capacityBeds}
& \sum_{p \in P_B} \sum_{q = \max\{1, h+1-f_p\}}^h y_{ptq}
   + \sum_{p \in P \setminus P_B} \sum_{q = \max\{1, h+1-f_p\}}^h z^B_{ptq}
\leq |B| 
& \forall t \in T, h \in H.
\end{align}
The objective function~\eqref{obf:max_p_surr} maximizes the number of treated patients, namely the patients assigned to a day and starting time. 
Constraints~\eqref{Cs:visit_assignment} guarantee that each patient is visited at most once in the planning horizon and are equivalent to~\eqref{C:visit_assignment}. 
Constraints~\eqref{Cs:onePerRoom} guarantee that the number of patients visited in each time slot is not greater than the number of exam rooms available. 
Constraints~\eqref{Cs:onePerRoom} are an aggregate version of \eqref{C:onePerRoom}. 
Constraints~\eqref{Cs:visit_infusion_1} and \eqref{Cs:visit_infusion_0} guarantee that a patient who is visited in a day receives the infusion in the same day,  for critical and non-critical patients, respectively. Constraints \eqref{Cs:first_visit_then_infusion_1} and \eqref{Cs:first_visit_then_infusion_0} guarantee that the infusion starts only after the consultation visit is finished. 
Constraints~\eqref{Cs:capacitychairs} and~\eqref{Cs:capacityBeds} guarantee that the number of patients who receive the infusion in the same time slot is not greater than the number of available chairs $|S|$ and beds $|B|$, respectively. 
They are aggregate versions of \eqref{C:capacitychairs} and \eqref{C:capacityBeds}, respectively.

We now prove that the aggregate formulation $(AF_1)$ is equivalent to the formulation \eqref{P1}. 
First, it is easy to show that from any feasible solution of \eqref{P1} it is possible to derive a feasible solution of $(AF_1)$ with the same objective function value.

\begin{lemma}
\label{lem:CF->SF}
If $(\alpha,\beta,\gamma^B,\gamma^S)$ is a feasible solution of \eqref{P1}, i.e., it satisfies constraints \eqref{C:visit_assignment}--\eqref{C:capacityBeds}, then the vector $(x,y,z^B,z^S)$ defined as follows:
\begin{align}
  & x_{pth} := \sum_{r \in R} \alpha_{pthr}, 
  & \forall\ p \in P, \ t \in T, \ h \in H_V,
  \label{e:x-alpha}
  \\
  & y_{pth} := \sum_{b \in B} \beta_{pthb},
  & \forall\ p \in P_B, \ t \in T, \ h \in H,
  \label{e:y-beta}
  \\
  & z^B_{pth} := \sum_{b \in B} \gamma^B_{pthb},
  & \forall\ p \in P \setminus P_B, \ t \in T, \ h \in H,
  \label{e:zB-gammaB}
  \\
  & z^S_{pth} := \sum_{s \in S} \gamma^S_{pths},
  & \forall\ p \in P \setminus P_B, \ t \in T, \ h \in H,
  \label{e:zS-gammaS}
\end{align}
is a feasible solution of $(AF_1)$, i.e., it satisfies constraints \eqref{Cs:visit_assignment}--\eqref{Cs:capacityBeds}, and the objective function values coincide: $F_1(x)=\varphi_1(\alpha)$.
\end{lemma}

\begin{proof}
It follows directly from constraints~\eqref{C:visit_assignment}  that $(x,y,z^B,z^S)$ satisfies constraints~\eqref{Cs:visit_assignment}.
The definition of variables $\alpha$ and constraints~\eqref{C:onePerRoom} guarantee that constraints~\eqref{Cs:onePerRoom} hold.
Constraints~\eqref{Cs:visit_infusion_1}--\eqref{Cs:first_visit_then_infusion_0} follow from constraints~\eqref{C:visit_infusion_1}--\eqref{C:first_visit_then_infusion_0}. 
By summing over $s \in S$ and $b \in B$ constraints~\eqref{C:capacitychairs} and \eqref{C:capacityBeds} we get constraints~\eqref{Cs:capacitychairs} and \eqref{Cs:capacityBeds}, respectively.
Finally, by definition of $x$, we get
$$
    F_1(x) = 
    \sum_{p \in P} \sum_{t \in T} \sum_{h \in H_V} x_{pth}
    =
    \sum_{p \in P} \sum_{t \in T} \sum_{h \in H_V} \sum_{r \in R} \alpha_{pthr}
    =
    \varphi_1(\alpha).
\vspace*{-7mm}
$$
\end{proof}
We now prove that a feasible solution of \eqref{P1}, which assigns a room and a bed (or a chair) to each patient, can be derived from a feasible solution of the aggregate formulation $(AF_1)$. 
This can be achieved by applying Algorithms~\ref{alg:1}, \ref{alg:2} and \ref{alg:3} described below.

\begin{algorithm}[htbp]
\small
	\DontPrintSemicolon
	\SetNoFillComment
	\caption{Assign patients to exam rooms}
	\label{alg:1}
	\KwData{vector $x=(x_{pth})$}
	\KwResult{vector $\alpha=(\alpha_{pthr})$}
    $R_{kt} \gets \{r \in R:\ w_{rkt}=1\}$ for any $k \in K$ and $t \in T$ 
    \tcp*{set of exam rooms for pathology $k$ in day $t$}
	$\alpha_{pthr} \gets 0$ for any $p \in P, t \in T, h \in H_V, r \in R$
	\tcp*{initialization}
	\For{$t \in T$}{		
		$L_r \gets 0$ for any $r \in R$ \tcp*{$L_r=1$ if room $r$ is occupied, $L_r=0$ if $r$ is empty} 
		\For{$h \in H_V$}{
			\tcc{Update the set of empty exam rooms}
			\For{$p \in P$}{
				\If{$x_{p,t,h-v_p}=1$}{
                \tcc{patient $p$ ends visit in slot $h$}
					\For{$r \in R$}{
						\If{$\alpha_{p,t,h-v_p,r}=1$}{
							$L_r \gets 0$ \tcp*{patient $p$ releases the room $r$}
						} 
					}
				}
			}
			\tcc{Assign patients to exam rooms}
			\For{$p \in P$}{
				\If{$x_{pth}=1$}{
					\For{$k \in K$}{
						\If{$u_{pk}=1$}{
                            \tcc{Assign $p$ to an empty room devoted to the corresponding pathology}
							$r \gets$ any element of the set $\{ r \in R_{kt}:\ L_r=0\}$ \;
							$\alpha_{pthr} \gets 1$ \;
							$L_r \gets 1$ \;
						}
					}
				}
			}
		}
	}
\end{algorithm}

Specifically, given the variables $x_{pth}$ satisfying the constraints of $(AF_1)$, Algorithm~\ref{alg:1} allows to define variables $\alpha_{pthr}$ satisfying the constraints of \eqref{P1}. 
In Algorithm~\ref{alg:1} we define $R_{kt}$ as the set of the exam rooms assigned to pathology $k$ in day $t$ (line 1) and introduce a binary label $L_r$ for each exam room $r$: $L_r=1$ if room $r$ is occupied at some time slot, and  $L_r=0$ if it is empty. 
All the variables $\alpha_{pthr}$ are initialized to 0 (line 2). 
For each day $t \in T$ the following steps are performed. 
First, all the labels $L_r$ are set to 0 (line 4) because all the exam rooms are empty at the beginning of the day.
For each time slot $h \in H_V$, first the set of empty exam rooms is updated: if a patient $p$ started the visit at slot $h-v_p$ (line 7) and was assigned to room $r$ (line 9), then $p$ releases the room $r$ at slot $h$, and thus its label $L_r$ is set to 0. 
Then, each patient who starts the visit at slot $h$ is assigned to an empty exam room devoted to the corresponding pathology (lines 15--25).

Similarly to Algorithm~\ref{alg:1}, Algorithm~\ref{alg:2} assigns non-critical patients to chairs (setting the variables $\gamma^S_{pths}$) starting from the variables $z^S_{pth}$ satisfying the constraints of the aggregate formulation $(AF_1)$.
In this case, a binary label $L_s$ is defined for each chair $s$: $L_s=1$ if chair $s$ is occupied at some time slot, and $L_s=0$ if it is empty. 
In the beginning, all the variables $\gamma^S_{pths}$ are initialized to 0 (line 1). 
For each day $t \in T$, all the labels $L_s$ are set to 0 (line 3) and, for each time slot $h \in H_V$, first the set of empty chairs is updated by checking which are the chairs released by patients who end the infusion at slot $h$ (lines 5--13). 
Then, each non-critical patient who starts the infusion at slot $h$ is assigned to an empty chair (lines 14--20).

\begin{algorithm}[htbp]
\small
	\DontPrintSemicolon
	\SetNoFillComment
	\caption{Assign patients to chairs}
	\label{alg:2}
	\KwData{vector $z^S=(z^S_{pth})$}
	\KwResult{vector $\gamma^S=(\gamma^S_{pths})$}

    $\gamma^S_{pths} \gets 0$ for any $p \in P, t \in T, h \in H, s \in S$
	\tcp*{initialization}
	\For{$t \in T$}{		
	$L_s \gets 0$ for any $s \in S$ \tcp*{$L_s=1$ if chair $s$ is occupied, $L_s=0$ if $s$ is empty} 
	\For{$h \in H$}{
		\tcc{Update the set of empty chairs}
		\For{$p \in P \setminus P_B$}{
			\If{$z^S_{p,t,h-f_p}=1$}{
   \tcc{patient $p$ ends infusion in slot $h$}
				\For{$s \in S$}{
					\If{$\gamma^S_{p,t,h-f_p,s}=1$}{
						$L_s \gets 0$ 
                        \tcp*{patient $p$ releases the chair $s$}
					}
				}
			}
		}
		\tcc{Assign patients to chairs}
		\For{$p \in P \setminus P_B$}{
			\If{$z^S_{pth}=1$}{
                        $s \gets$ any element of the set $\{ s \in S:\ L_s=0\}$ \;
						$\gamma^S_{pths} \gets 1$
                        \tcp*{Assign patient $p$ to an empty chair}
						$L_s \gets 1$ \;
				}
			}

		}
	}
\end{algorithm}

Finally, Algorithm~\ref{alg:3} assigns critical and non-critical patients to beds, setting the variables $\beta_{pthb}$ and $\gamma^B_{pthb}$, respectively, starting from variables $y_{pth}$ and $z^B_{pth}$ satisfying the constraints of formulation $(AF_1)$.
The steps of Algorithm~\ref{alg:3} are similar to those of Algorithms~\ref{alg:1} and~\ref{alg:2}: for each day $t$ and each time slot $h$, we first update the set of empty beds, and then we assign patients to beds.

\begin{algorithm}[htbp]
\small
	\DontPrintSemicolon
	\SetNoFillComment
	\caption{Assign patients to beds}
	\label{alg:3}
	\KwData{vectors $y=(y_{pth})$ and $z^B=(z^B_{pth})$}
	\KwResult{vectors $\beta=(\beta_{pthb})$ and $\gamma^B=(\gamma^B_{pthb})$}

    $\beta_{pthb} \gets 0$ for any $p \in P_B, t \in T, h \in H, b \in B$
	\tcp*{initialization}
    $\gamma^B_{pthb} \gets 0$ for any $p \in P \setminus P_B, t \in T, h \in H, b \in B$
    
	\For{$t \in T$}{		
		$L_b \gets 0$ for any $b \in B$ \tcp*{$L_b=1$ if bed $b$ is occupied, $L_b=0$ if $b$ is empty} 
		\For{$h \in H$}{
		\tcc{Update the set of empty beds}
\For{$p \in P$}{
	\If{$p \in P_B$ {\bf and} $y_{p,t,h-f_p}=1$}{
    \tcc{critical patient $p$ ends infusion in slot $h$}
		\For{$b \in B$}{
			\If{$\beta_{p,t,h-f_p,b}=1$}{
				$L_b \gets 0$ 
                \tcp*{critical patient $p$ releases the bed $b$}
			}
		}
	}
	\If{$p \in P \setminus P_B$ {\bf and} $z^B_{p,t,h-f_p}=1$}{
    \tcc{non-critical patient $p$ ends infusion in slot $h$}
		\For{$b \in B$}{
			\If{$\gamma^B_{p,t,h-f_p,b}=1$}{
				$L_b \gets 0$
                \tcp*{non-critical patient $p$ releases the bed $b$}
			}
		}
	}
}
\tcc{Assign patients to beds}
\For{$p \in P$}{
	\If{$p \in P_B$ {\bf and} $y_{pth}=1$}{
		$b \gets$ any element of the set $\{ b \in B:\ L_b=0\}$ \;
		$\beta_{pthb} \gets 1$
        \tcp*{Assign critical patient $p$ to an empty bed}
		$L_b \gets 1$ \;
	}
	\If{$p \in P \setminus P_B$ {\bf and} $z^B_{pth}=1$}{
		$b \gets$ any element of the set $\{ b \in B:\ L_b=0\}$ \;
		$\gamma^B_{pthb} \gets 1$
        \tcp*{Assign non-critical patient $p$ to an empty bed}
		$L_b \gets 1$ \;
	}
}
		}
	}
\end{algorithm}

\begin{lemma}
\label{lem:SF->CF}
If $(x,y,z^B,z^S)$ is a feasible solution of $(AF_1)$, i.e., it satisfies constraints \eqref{Cs:visit_assignment}--\eqref{Cs:capacityBeds}, and we apply Algorithms~\ref{alg:1}, \ref{alg:2} and \ref{alg:3}, then we find a vector $(\alpha,\beta,\gamma^B,\gamma^S)$ that is feasible for \eqref{P1}, i.e., it satisfies constraints \eqref{C:visit_assignment}--\eqref{C:capacityBeds}, and $\varphi_1(\alpha)=F_1(x)$.
\end{lemma}

\begin{proof}
First, we notice that the vector $(\alpha,\beta,\gamma^B,\gamma^S)$ found applying Algorithms~\ref{alg:1}, \ref{alg:2} and \ref{alg:3} satisfies relations 
\eqref{e:x-alpha}--\eqref{e:zS-gammaS} with the vector $(x,y,z^B,z^S)$. 
In fact, if for some patient $p \in P$ the variables $x_{pth}=0$ for any $t \in T$ and $h \in H_V$, then all the corresponding variables $\alpha_{pthr}=0$ for any $t \in T$, $h \in H_V$ and $r \in R$ since the condition at line 16 of Algorithm~\ref{alg:1} is never satisfied. 
Otherwise, if $x_{pth}=1$ for some $t \in T$ and $h \in H_V$, then there exists a unique $r \in R$ such that $\alpha_{pthr}=1$, while $\alpha_{pthr'}=0$ for any $r' \in R \setminus\{r\}$ (line 19 of Algorithm~\ref{alg:1}). Hence, relation~\eqref{e:x-alpha} holds. 
Similarly, relations~\eqref{e:y-beta}--\eqref{e:zS-gammaS} can be proved. 
Therefore, $(\alpha,\beta,\gamma^B,\gamma^S)$ satisfies constraints~\eqref{C:visit_assignment} and \eqref{C:visit_infusion_1}--\eqref{C:first_visit_then_infusion_0} since $(x,y,z^B,z^S)$ is a feasible solution of $(AF_1)$.
Moreover, the use of labels $L_r$ in Algorithm~\ref{alg:1}, to denote which rooms are occupied during a day, allows to have at most one patient who is visited in each time slot in each room. Thus, constraints~\eqref{C:onePerRoom} are satisfied as well.
Similarly, the use of labels $L_s$ in Algorithm~\ref{alg:2} and $L_b$ in Algorithm~\ref{alg:3} allow to have at most one patient who receives the infusion in each time slot and in each chair or bed, respectively. 
Therefore, also constraints~\eqref{C:capacitychairs} and \eqref{C:capacityBeds} are satisfied. 
Thus, $(\alpha,\beta,\gamma^B,\gamma^S)$ is a feasible solution of the formulation \eqref{P1}. 
Finally, relation~\eqref{e:x-alpha} implies that $\varphi_1(\alpha)=F_1(x)$ holds.
\end{proof}

The equivalence between the complete and aggregate formulations  follows from Lemma~\ref{lem:CF->SF} and Lemma~\ref{lem:SF->CF}.

\begin{theorem}
\label{t:CF=SF}
The aggregate formulation $(AF_1)$ is equivalent to formulation~\eqref{P1}.
\end{theorem}


\subsection{Solving problem $(P_2)$}
\label{s:solvingP2}

Similarly to the aggregate formulation of problem $(P_1)$, we consider the aggregate formulation of problem $(P_2)$, which is defined as follows:
\begin{align}
\label{SF2}
\tag{$AF_2$} 
\left\{
\begin{array}{ll}
\min\limits_{(x,y,z^B,z^S)}\ & F_2(x,y,z^B,z^S) 
\\
\text{subject to} & \eqref{Cs:visit_assignment}-\eqref{Cs:capacityBeds}
\\
& F_1(x) \geq \bar{v}_1
\end{array}
\right.
\end{align}
where $\bar{v}_1$ is the optimal value of $(AF_1)$, and
\begin{equation}
\label{OBJ:waitingTime-xyz}
F_2(x,y,z^B,z^S) = 
\sum_{t \in T}
F^t_2(x,y,z^B,z^S)
\end{equation}
where
\begin{equation}
\label{e:F2t}
\begin{array}{rl}
F^t_2(x,y,z^B,z^S) = 
\max & \left\{ 
\displaystyle
\max_{p \in P_B} \left\{ 
 \sum_{h \in H} h y_{pth} 
- 
\sum_{h \in H_V} (h+v_p) x_{pth}
\right\}
\right.
, 
\\[4mm]
& 
\left.
\max\limits_{p \in P \setminus P_B} \left\{ 
\sum\limits_{h \in H} h \left[
z^B_{pth} + z^S_{pth}
\right]
- 
\sum\limits_{h \in H_V} (h+v_p) x_{pth}
\right\}
\right\}
\end{array}
\end{equation}
represents the maximum waiting time in day $t$.
Following the same reasoning made in Section~\ref{s:solvingP1}, it is possible to prove that the aggregate formulation $(AF_2)$ is equivalent to the formulation \eqref{P2} of problem $(P_2)$.

However, the aggregate formulation $(AF_2)$ proves to be more computationally challenging than $(AF_1)$ (see Section~\ref{sec::results-P2}). 
Therefore, we developed a procedure to solve $(AF_2)$, which is reported in Procedure~\ref{proc:1}.


\begin{procedure}[tbp]
\small
	\DontPrintSemicolon
	\SetNoFillComment
	\caption{1(): Finding an optimal solution of the formulation $(AF_2)$}
	\label{proc:1}
    Find an optimal solution $(\bar{x},\bar{y},\bar{z}^B,\bar{z}^S)$ of the formulation $(AF_1)$ \;
    $\bar{v}_1 \gets F_1(\bar{x})$ \;
    \tcc{Fix the set of treated patients each day}
    \For{$t \in T$}{	
        $P_t \gets \left\{p \in P:\ \sum\limits_{h \in H_V} \bar{x}_{pth} = 1 \right\}$ \;
    }
    \tcc{Once the set of treated patients and their assigned day are fixed, solve formulation $(AF_2)$ by decomposing it into a set of problems one for each day}
    \For{$\tau \in T$}{	
        find an optimal solution of the problem:
        $$
        \left\{
        \begin{array}{ll}
        \min\limits_{(x,y,z^B,z^S)}\ & F_2(x,y,z^B,z^S) 
        \\
        \text{subject to} & \eqref{Cs:visit_assignment}-\eqref{Cs:capacityBeds}
        \\
                    & F_1(x) \geq |P_{\tau}|
        \end{array}
        \right.
        $$
        where the variables $x_{pth}$, $y_{pth}$, $z^B_{pth}$ and $z^S_{pth}$ are restricted to $p \in P_{\tau}$, and $t=\tau$ \; 
    }
   \tcc{Solve formulation $(AF_2)$ with warm start}
    Use the union of the solutions found at lines 6-8 as a warm start to solve the whole formulation $(AF_2)$:
        $$
        \left\{
        \begin{array}{ll}
        \min\limits_{(x,y,z^B,z^S)}\ & F_2(x,y,z^B,z^S) 
        \\
        \text{subject to} & \eqref{Cs:visit_assignment}-\eqref{Cs:capacityBeds}
        \\
                    & F_1(x) \geq \bar{v}_1
        \end{array}
        \right.
        $$

\end{procedure}

First, we solve aggregate formulation $(AF_1)$ of problem $(P_1)$ (line 1) and fix the set of treated patients and their assigned day (lines 3-5). 
Then, we solve the aggregate formulation $(AF_2)$ of problem $(P_2)$, where the day of treatment of each treated patient is set. 
This latter problem can be decomposed into a set of $|T|$ problems, one for each day, which optimize the starting time of visits and infusions, with the goal of minimizing the maximum waiting time of patients (lines 6-8). 

Finally, the obtained solution is used as a warm start for solving the whole formulation~\eqref{SF2} (line 9). 
Notice that in the last step, we impose that the total number of treated patients is optimal ($F_1(x) \geq \bar{v}_1$), but we allow to change (with respect to the optimal solution of $(AF_1)$) the set of treated patients as well as their visit and infusion days and times.


\subsection{Solving problem $(P_3)$}
\label{s:solvingP3}

Problem $(P_3)$ proves even more challenging than $(P_2)$. 
We devised a two steps procedure that computes a heuristic solution applying a decomposition and a $k$-opt neighborhood search. Further, we derived two upper bounds for the problem as a comparison for the heuristic solution.


\subsubsection{Heuristic for $(P_3)$}

Similarly to the aggregate formulation of problems $(P_1)$ and $(P_2)$, we consider the aggregate formulation of problem $(P_3)$, which is defined as follows:
\begin{align}
\label{SF3}
\tag{$AF_3$} 
\left\{
\begin{array}{lll}
\max\limits_{(x,y,z^B,z^S)}\ & F_3(z^S) &  
\\
\text{subject to} & \eqref{Cs:visit_assignment}-\eqref{Cs:capacityBeds} &
\\
& F_1(x) \geq \bar{v}_1 &
\\
& F^t_2(x,y,z^B,z^S) \leq \bar{v}^t_2 & \qquad \forall\ t \in T,
\end{array}
\right.
\end{align}
where $\bar{v}_1$ is the optimal value of $(AF_1)$, $\bar{v}^t_2$ is the maximum waiting time in day $t$ given by the optimal solution of $(AF_2)$ and
\begin{align} 
\label{OBJ:chairPreferences-xyz}
    F_3(z^S) = 
    \sum_{p \in P \setminus P_B} \sum_{t \in T} \sum_{h \in H} z^S_{pth}.
\end{align}
Following the same reasoning made in Section~\ref{s:solvingP1}, it is possible to prove that the aggregate formulation $(AF_3)$ is equivalent to formulation~\eqref{P3}.

\begin{procedure}[htbp]
\small
\DontPrintSemicolon
	\SetNoFillComment
	\caption{2(): Finding a solution of the  formulation $(AF_3)$}
	\label{proc:2}
    Apply Procedure~\ref{proc:1} to find an optimal solution $(\bar{x},\bar{y},\bar{z}^B,\bar{z}^S)$ of the formulation $(AF_2)$ \;
    \tcc{Fix the maximum waiting time and the set of treated patients in each day}
    \For{$t \in T$}{
        $\bar{v}^t_2 \gets F^t_2(\bar{x},\bar{y},\bar{z}^B,\bar{z}^S)$
        \tcp*{maximum waiting time in day $t$, see \eqref{e:F2t}}
        $P_t \gets \left\{p \in P:\ \sum\limits_{h \in H_V} \bar{x}_{pth} = 1 \right\}$
        \tcp*{set of treated patients in day $t$}
    }
    \tcc{Once the set of treated patients and their assigned day are fixed, solve the aggregate formulation $(AF_3)$ by decomposing it into a set of problems one for each day}
    \For{$\tau \in T$}{	
        find an optimal solution of the single day problem:
        $$
        \left\{
        \begin{array}{ll}
        \max\limits_{(x,y,z^B,z^S)}\ & F_3(z^S) 
        \\
        \text{subject to} & \eqref{Cs:visit_assignment}-\eqref{Cs:capacityBeds}
        \\
            & F_1(x) \geq |P_{\tau}|
        \\
            & F^{\tau}_2(x,y,z^B,z^S) \leq \bar{v}^{\tau}_2    
        \end{array}
        \right.
        $$
        where the variables $x_{pth}$, $y_{pth}$, $z^B_{pth}$ and $z^S_{pth}$ are restricted to $p \in P_{\tau}$, and $t=\tau$\; 
    }
    \tcc{$k$-opt neighborhood search}
    Use the union of the solutions found at lines 6-8 as starting solution of a $k$-opt neighborhood search for the whole formulation $(AF_3)$:
    $$ 
\left\{
\begin{array}{lll}
\max\limits_{(x,y,z^B,z^S)}\ & F_3(z^S) &  
\\
\text{subject to} & \eqref{Cs:visit_assignment}-\eqref{Cs:capacityBeds} &
\\
& F_1(x) \geq \bar{v}_1 &
\\
& F^t_2(x,y,z^B,z^S) \leq \bar{v}^t_2 & \qquad \forall\ t \in T
\end{array}
\right.
$$
\end{procedure}
We find a feasible solution of $(AF_3)$ by means of Procedure~\ref{proc:2}.
First, we apply Procedure~\ref{proc:1} to solve $(AF_2)$ (line 1) and fix the daily maximum waiting time, the set of treated patients, and their assigned day (lines 2-5). 
Then, we solve \eqref{SF3}, where the day of treatment of each treated patient is set. 
This latter problem can be decomposed into a set of $|T|$ problems, one for each day, which optimize the starting time of the visits and infusions, with the goal of maximizing the number of non-critical patients who receive the treatment on a chair (lines 6-8). 
Finally, the obtained solution is used as a starting solution for a $k$-opt neighborhood search applied to the whole aggregate formulation $(AF_3)$ (line 9).

In the $k$-opt neighborhood search we consider, in each iteration, a neighborhood of the current solution that is defined as the set of solutions that differ from the current one by the value of a limited number of variables~\citep{FL2003}. Rather than explicitly enumerating all the neighbors, the neighborhood is searched solving a MILP model. 
The following $k$-opt neighborhood constraints are added to limit the size of the neighborhood, namely the difference between the current solution and the neighbor ones:
\begin{align}
\label{con-p3-koptX}
& \sum_{p \in P} \sum_{t \in T} \sum_{\substack{h \in H: \\ \bar{x}_{pth}=1}}(1-x_{pth}) +
\sum_{p \in P} \sum_{t \in T} \sum_{\substack{h \in H: \\ \bar{x}_{pth}=0}} x_{pth} \leq k_x,
\\
\label{con-p3-koptY}
& \sum_{p \in P_B} \sum_{t \in T} \sum_{\substack{h \in H: \\ \bar{y}_{pth}=1}}(1-y_{pth}) +
\sum_{p \in P_B} \sum_{t \in T} \sum_{\substack{h \in H: \\ \bar{y}_{pth}=0}} y_{pth} \leq k_y,
\\
\label{con-p3-koptZB}
& \sum_{p \in P\setminus P_B} \sum_{t \in T} \sum_{\substack{h \in H: \\ \bar{z}^B_{pth}=1}}(1-z^B_{pth}) +
\sum_{p \in P\setminus P_B} \sum_{t \in T} \sum_{\substack{h \in H: \\ \bar{z}^B_{pth}=0}} z^B_{pth} \leq k_z^B,
\\
\label{con-p3-koptZS}
& \sum_{p \in P\setminus P_B} \sum_{t \in T} \sum_{\substack{h \in H: \\ \bar{z}^S_{pth}=1}}(1-z^S_{pth}) +
\sum_{p \in P\setminus P_B} \sum_{t \in T} \sum_{\substack{h \in H: \\ \bar{z}^S_{pth}=0}} z^S_{pth} \leq k_z^S,
\end{align}
where $\bar{x}_{pth}, \bar{y}_{pth}, \bar{z}^B_{pth}, \bar{z}^S_{pth}$ represent the current solution.
The size of the neighborhood is limited by parameters $k_x, k_y, k_z^B, k_z^S$, respectively.
In each iteration, we look for the best solution w.r.t. the objective function \eqref{OBJ:chairPreferences-xyz} in a steepest descent fashion. 
However, the neighborhood is searched heuristically, as a time limit is set for each MILP problem. The current solution is replaced by the solution provided by the MILP if the latter improves upon the former. The $k$-opt neighborhood step can modify all the variables, while keeping the values of $F_1$ and $F^t_2$ bounded by $\bar{v}_1$ and $\bar{v}^t_2$, respectively. 


\subsubsection{Upper bounds for $(P_3)$}\label{sec::upperBounds}

As mentioned, to assess the quality of the heuristic solution for $(AF_3)$, we devised two upper bounds of the maximum number of patients who can receive the infusion in a chair.

\subsubsection*{Multiple knapsack based upper bound $(UB_1)$}
\label{sec:UB1}

The first upper bound is based on a multiple knapsack problem, where each chair represents a knapsack and patients represent the items. 
We assume that each patient $p$ has a weight equal to $f_p$ and a profit equal to 1. 
The capacity of each chair is the number of slots that can be used for infusion. Although the chairs are available for $|H|$ time slots, in the first time slots of each day no patients can use them, as no patient has undergone the medical consultation, yet. 
So, on each day, all the chairs are empty and idle in the first $\min\limits_{p \in P\setminus P_B} v_p$ time slots. Therefore, the capacity of each chair is set equal to $(|H|-\min\limits_{p \in P\setminus P_B} v_p)$. 

The multiknapsack based upper bound $(UB_1)$ can be formulated with the following binary variables:
$$
    \mu_{ps}=
    \begin{cases} 
    1 & \text{if patient $p \in P \setminus P_B$ is assigned to chair $s$,} 
    \\
    0 & \text{otherwise,}
    \end{cases}
$$
and modeled as follows:
\begin{align}
\label{UB1:fob}& \max\ \sum_{p \in P \setminus P_B} \sum_{s \in S} \mu_{ps} &
\\
\label{UB1:assignment}& \sum_{s \in S} \mu_{ps} \leq 1 
& \forall\ p \in P \setminus P_B,
\\
\label{UB1:capacity}& \sum_{p \in P \setminus P_B} f_p \, \mu_{ps} \leq |T| \, \left(
|H| - \min_{p \in P \setminus P_B} v_p
\right)
& \forall\ s \in S.
\end{align}
The objective function~\eqref{UB1:fob} aims at maximizing the number of patients assigned to the chairs while guaranteeing that the chair capacities are not exceeded~\eqref{UB1:capacity}. Constraints~\eqref{UB1:assignment} guarantee that each patient is assigned to at most one chair. 


\subsubsection*{Improved upper bound $(UB_2)$}\label{sec:UB2}

We now formulate an alternative problem to compute an improved upper bound of the number of patients who can receive the infusion in a chair, by considering additional features of the original problem to improve the bound w.r.t. $(UB_1)$. 
In problem $(UB_2)$, we assign patients to single chairs and beds every day, similarly to formulations~\eqref{P1}--\eqref{P3}, but we do not consider the visits and the scheduling decisions. 
Accounting for the multi-objective structure, the total number of visited patients is given by $\bar{v}_1$, i.e., the maximum number of visited patients computed by $(AF_1)$. 
The availability of skilled clinicians to treat the different pathologies is taken into account, as well as the assignment of pathologies to working days.
 
We model the assignment choices in problem $(UB_2)$ with the following three families of binary variables:
\begin{align*}
    & \lambda_{ptb}=
    \begin{cases} 
    1 & \text{if patient $p \in P_B$ is assigned to bed $b$ in day $t$},
    \\
    0 & \text{otherwise},
    \end{cases}
\\
    & \mu^B_{ptb}=
    \begin{cases} 
    1 & \text{if patient $p \in P \setminus P_B$ is assigned to bed $b$ in day $t$},
    \\
    0 & \text{otherwise},
    \end{cases}
\\
    & \mu^S_{pts}=
    \begin{cases} 
    1 & \text{if patient $p \in P \setminus P_B$ is assigned to chair $s$ in day $t$,} 
    \\
    0 & \text{otherwise.}
    \end{cases}
\end{align*}
The variables are defined so as to prevent the assignment of critical and non-critical patients to beds and chairs in days where there are no exam rooms devoted to the patient's pathology (i.e., $\sum_{r \in R} w_{rkt}=0$). 

The chair capacities can be further refined w.r.t. the capacity of the multiknapsack problem, namely $\left( |H|-\min\limits_{p \in P\setminus P_B} v_p \right)$, both at the beginning and at the end of the day. 

\begin{figure}[tb]
\begin{center}
\includegraphics[width=0.5\textwidth]{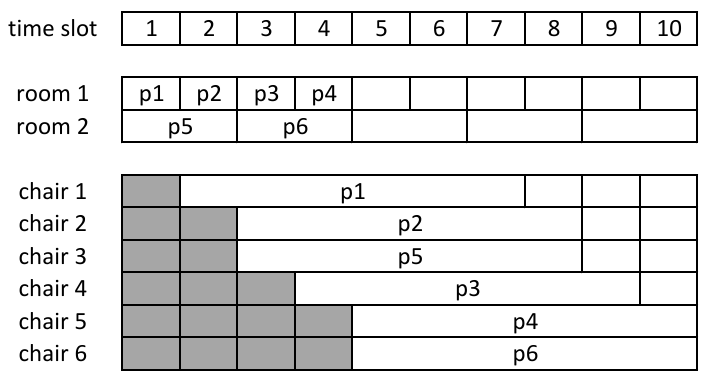}
\caption{Unavailable slots in each chair at the beginning of the day.
}
\label{fig:1}
\end{center}
\end{figure}

Let us consider first the beginning of the day. 
Figure~\ref{fig:1} represents an example with 10 time slots, 6 patients, 2 rooms, and 6 chairs. 
Let us suppose that the visit time of the first four patients (in room 1) is equal to one slot, while the visit time of the last two patients (in room 2) is equal to two slots. 
All the patients have the same infusion time, equal to 6 slots. 
As discussed for $(UB_1)$, in the first time slot of the day, all the chairs are empty. 
However, as the exam rooms are only two, at most two patients can start their infusion in the second time slot, thus the remaining four chairs cannot be used in the second time slot, either (actually, only patient p1 can start the infusion in the second time slot).
Similarly, in the third time slot, at most four chairs can be used, as  at most two additional patients have started and finished their visits and can start the infusion. 
By generalizing the reasoning, it is possible to better evaluate the actual capacity of each chair. 
Chairs are divided into $\left\lceil \frac{|S|}{|R|} \right\rceil$ subsets: each subset $j = 1,\dots,\left\lceil \frac{|S|}{|R|} \right\rceil$ contains the chairs that can actually be used in the same set of time slots. 
Therefore, the capacity of each chair belonging to the subset $j$ is reduced by $j \left( \min\limits_{p \in P \setminus P_B} v_p \right)$ slots. 
For subset $1$, which contains $|R|$ chairs, all the time slots but $\min\limits_{p \in P \setminus P_B} v_p$ are available and can be used; for the subset $2$, which contains $|R|$ chairs, all the time slots but $2\left(\min\limits_{p \in P \setminus P_B} v_p\right)$ can be actually used, and so on.
Thus, constraints~\eqref{UB1:capacity} can be replaced by
\begin{align*}
\begin{array}{l}
\sum\limits_{p \in P \setminus P_B} f_p \mu^S_{pts} 
\leq 
|H| - j \left( \min\limits_{p \in P \setminus P_B} v_p \right) \qquad \qquad
\end{array}
& 
\begin{array}{l}
\forall\ t \in T, \
j = 1,\dots,\left\lceil \frac{|S|}{|R|} \right\rceil,
\\
s = 1+(j-1)|R|,\dots,\min\{j|R|,\, |S|\}.
\end{array}
\end{align*}
The capacity of each chair can be further reduced at the end of the day taking into account the impact of the maximum waiting time in day $t$ as set by problem $(P_2)$, which limits the initial time of infusions, which can start at last at time slot $|H_V| + \bar{v}^t_2$ in day $t$. 

\begin{figure}[tbp]
\begin{center}
\includegraphics[width=\textwidth]{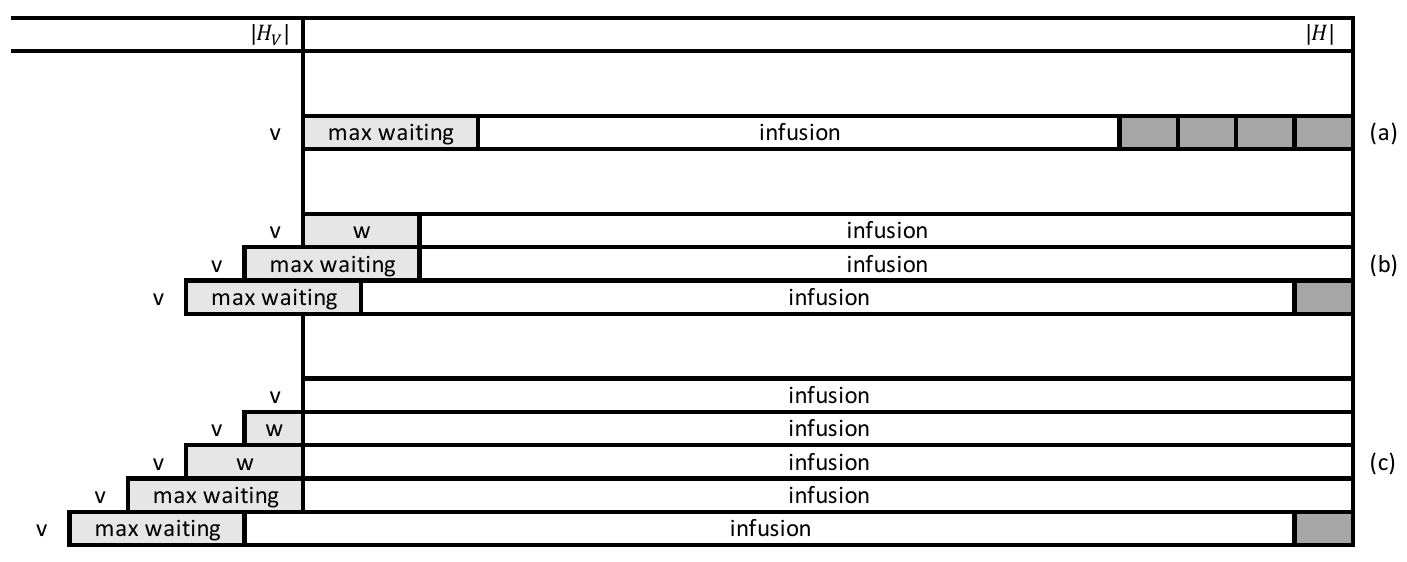}
\caption{Unused slots in each chair at the end of the day.}
\label{fig:2}
\end{center}
\end{figure}

Let us consider Figure~\ref{fig:2}, which represents the three cases that may occur, depending on the length of infusion $f_p$.

In case (a) we consider the patients whose infusion time is such that $f_p < |H| - (|H_V| + \bar{v}^t_2)$. For such patients, the infusion ends before $|H|$ even if it starts in the last possible time slot, and therefore the chair remains empty in $|H| - (|H_V| + \bar{v}^t_2 + f_p)$ time slots at the end of the day (grey slots).

In case (b) we consider the patients whose infusion time is such that 
$$
    |H| - (|H_V| + \bar{v}^t_2) \leq f_p < |H| - |H_V|.
$$
Their infusions can fully use the time slots of the day provided that they wait for $0< w \leq \bar{v}^t_2$ slots between the end of the visit and the beginning of the infusion. However, due to the limited number of rooms, not all such patients can start their visit so late as to end precisely in the last time slot $|H|$. Instead, some of them may have to start before time slot $|H| - ( f_p + \bar{v}^t_2)$, thus leaving some of the time slots unused at the end of the day (last patient in part (b) of the Figure).

In case (c) we consider the patients whose infusion time is such that $f_p \geq |H| - |H_V|$. Those patients, even if they start their infusion right after the end of the visit, can fully exploit the time slots at the end of the day provided that they do not end their visit too early. Note that, due to the limited number of exam rooms, the number of such patients may not be sufficient to fully use the last time slots of the day. For example, in Figure~\ref{fig:2}, the first 4 patients end their infusion at the end of the day, while the fifth one starts and ends the visit too early and, therefore, even considering the maximum waiting time, leaves one time slot empty at the end of the day.

\begin{algorithm}[tbp]
\small
	\DontPrintSemicolon
	\SetNoFillComment
	\caption{Computing parameters $N^S_{it}$}
	\label{alg:6}
    \tcc{Compute the number $L^S_{\ell r t}$ of non-critical patients whose infusion time is equal to $\ell$ and can be visited in day $t$ in room $r$}
    \For{$\ell = 1,\dots,\max\limits_{p \in P \setminus P_B} f_p$, $r \in R$, $t \in T$}{
                $L^S_{\ell rt} \gets \left| \left\{ p \in P \setminus P_B:\ f_p=\ell \text{ and } \exists\ k \in K:\ u_{pk}=1,\ w_{rkt}=1 \right\} \right|$
    }
    \tcc{Compute the number $M^S_{\ell rt}(i)$ of patients who can leave $i$ time slots empty if they are visited as last, for each length of infusion $\ell$, room $r$, and day $t$}
    \For{$\ell = 1,\dots,\max\limits_{p \in P \setminus P_B} f_p$, $r \in R$, $t \in T$}{
                \uIf{$\ell \geq |H|-(|H_V|+\bar{v}^t_2)$}{
                    $M^S_{\ell rt}(0) \gets \min\left\{ L^S_{\ell rt} , \ 
                    \bar{v}^t_2 + 1 + \ell - (|H|-|H_V|) , \ \bar{v}^t_2 + 1
                    \right\}$
                    }
                \Else{$M^S_{\ell rt}(0) \gets 0$}
                \For{$i=1,\dots,|H|$}{ 
                    \uIf{$\ell \geq |H|-(|H_V|+ \bar{v}^t_2)-i$  {\bf and} $\sum\limits_{j=0}^{i-1} M^S_{\ell rt}(j) < L^S_{\ell rt}$}{
                        $M^S_{\ell rt}(i) \gets 1$
                    }
                    \Else{$M^S_{\ell rt}(i) \gets 0$}
                }
            }
    \tcc{Compute the number $N^S_{it}$ of chairs that are empty for the last $i$ time slots in day $t$}
    $N^S_{it} \gets 0$ for any $i=0,\dots,|H|$ and $t \in T$ \;
    \For{$t \in T$}{
        $\sigma \gets |S|$ \;
        $i \gets 0$ \;
        \While{$\sigma>0$ {\bf and} $i \leq |H|$}{
            $N^S_{it} \gets \min\left\{ \sigma , \ \sum\limits_{\ell = 1}^{\max\limits_{p \in P \setminus P_B} f_p} \sum\limits_{r \in R} M^S_{\ell rt}(i)  \right\}$ \;
            $\sigma \gets \sigma - N^S_{it}$ \;
            $i \gets i+1$ \;
        }
    }
\end{algorithm}
Based on the above cases, we can evaluate the number $N^S_{it}$ of chairs that are empty for the last $i$ time slots in day $t$, as described in Algorithm~\ref{alg:6}. 
As we aim at computing a relaxation, we fill as much as possible the last slots. 
The number $N^S_{0t}$ of chairs that are fully used  is the maximum possible, and the remaining patients are used to fill the other chairs' last slots.

First, lines 1-3 compute the number $L^S_{\ell r t}$ of non-critical patients whose infusion time is equal to $\ell$ and can be visited in day $t$ in room $r$,  
accounting for the availability of clinicians each day ($u_{pk}=1$) and the assignment of pathology to the rooms and days ($w_{rkt}=1$). 

Next, we compute the number $M^S_{\ell rt}(i)$ of patients who can leave $i$ time slots empty if they are visited as last, for each length of infusion $\ell$, room $r$ and day $t$ (Lines 4-17). 
For instance, $M^S_{\ell rt}(0)$ is the number of patients who use all the time slots if visited in the last slots available for visits. 
Parameter $M^S_{\ell rt}(0)$ is computed in lines 5-9, for a given room $r$. 
In particular, lines 5-6 refer to cases (b) and (c), while line 8 refers to case (a) in Figure~\ref{fig:2}. 
Lines 10-16 compute the values of $M^S_{\ell rt}(i)$ for $i \geq 1$. 

Finally, the value of $N^S_{it}$ is computed in lines 18-27, based on the $M^S_{\ell rt}(i)$ parameters.  
Similar reasoning can be applied for the beds, yielding to parameter $N^B_{it}$. 

Both for the chairs and the beds, we cannot decide which chairs (or beds) have a reduced capacity: instead, we let the optimal solution of problem $(UB_2)$ decide to which chairs (or beds) must be assigned the above computed reduced capacities. 
To formulate problem $(UB_2)$, we introduce the following binary variables:
$$
    \rho^B_{itb}=
    \begin{cases} 
    1 & \text{if bed $b$ has $i$ unused time slots at the end of day $t$},
    \\
    0 & \text{otherwise},
    \end{cases}
$$
$$
    \rho^S_{its}=
    \begin{cases} 
    1 & \text{if chair $s$ has $i$ unused time slots at the end of day $t$,} 
    \\
    0 & \text{otherwise.}
    \end{cases}
$$
The capacity of each chair $s$ is then $|H| - j \left( \min\limits_{p \in P \setminus P_B} v_p \right) 
- \sum\limits_{i=0}^{|H|} i \rho^S_{its}$.
The formulation of problem $(UB_2)$ is the following:
\begin{align}
\label{UB2:of}
& \max\ \sum_{p \in P \setminus P_B} \sum_{t \in T} \sum_{s \in S} \mu^S_{pts} &
\\
& \sum_{p \in P_B} \sum_{t \in T} \sum_{b \in B} \lambda_{ptb} + 
\sum_{p \in P \setminus P_B} \sum_{t \in T} 
\left[ 
\sum_{b \in B} \mu^B_{ptb} + \sum_{s \in S} \mu^S_{pts}
\right]
= \bar{v}_1
\label{UB2:v1}
& 
\\
& \sum_{t \in T} \sum_{b \in B} \lambda_{ptb} \leq 1 
& \forall p \in P_B 
\label{UB2:lambda}
\\
& \sum_{t \in T} \left[ 
\sum_{b \in B} \mu^B_{ptb} + \sum_{s \in S} \mu^S_{pts}
\right]
\leq 1 
& \forall p \in P \setminus P_B 
\label{UB2:sum_mu_1}
\end{align}

\begin{align}
& 
\begin{array}{l}
\sum\limits_{p \in P_B} f_p \lambda_{ptb} 
+ \sum\limits_{p \in P \setminus P_B} f_p \mu^B_{ptb} 
\\
\qquad \leq 
|H| - j \left( \min\limits_{p \in P} v_p \right) 
- \sum\limits_{i=0}^{|H|} i \rho^B_{itb}
\end{array}
& 
\begin{array}{l}
\forall t \in T, \
j = 1,\dots,\left\lceil \frac{|B|}{|R|} \right\rceil,
\\
b = 1+(j-1)|R|,\dots,\min\{j|R|,\, |B|\}
\label{UB2:bed_capacity}
\end{array}
\\
&
\sum_{i=0}^{|H|} \rho^B_{itb}=1
\label{UB2:roB=1}
&
\forall t \in T, \ b \in B
\\
&
\sum_{b \in B} \rho^B_{itb} = N^B_{it}
&
\forall i = 0,\dots,|H|, \ t \in T
\label{UB2:roB=Nb}
\\
& 
\begin{array}{l}
\sum\limits_{p \in P \setminus P_B} f_p \mu^S_{pts} 
\leq 
|H| - j \left( \min\limits_{p \in P \setminus P_B} v_p \right) 
- \sum\limits_{i=0}^{|H|} i \rho^S_{its}
\end{array}
& 
\begin{array}{l}
\forall t \in T, \
j = 1,\dots,\left\lceil \frac{|S|}{|R|} \right\rceil,
\\
s = 1+(j-1)|R|,\dots,\min\{j|R|,\, |S|\}
\label{UB2:chair_capacity}
\end{array}
\\
&
\sum_{i=0}^{|H|} \rho^S_{its} = 1
&
\forall t \in T, \ s \in S
\label{UB2:roS=1}
\\
&
\sum_{s \in S} \rho^S_{its} = N^S_{it}
&
\forall i = 0,\dots,|H|, \ t \in T.
\label{UB2:roS=Ns}
\end{align}
Objective function~\eqref{UB2:of} represents the number of patients assigned to chairs. Constraint~\eqref{UB2:v1} sets the total number of patients to assign to chairs and beds. Constraints~\eqref{UB2:lambda} ensure that each critical patient is assigned to at most one bed, while~\eqref{UB2:sum_mu_1} ensure that each non-critical patient is assigned at most once, either to a bed or to a chair.  
Constraints~\eqref{UB2:bed_capacity} and \eqref{UB2:chair_capacity} are the capacity constraints for each bed and for each chair, respectively, in terms of actual number of available time slots. Constraints~\eqref{UB2:roB=1} assign the number of unused slots to each bed, while~\eqref{UB2:roB=Nb} force the number of beds with $i$ empty slots to be equal to $N^B_{it}$, for each $i=1, \dots,|H|$, thus guaranteeing that the empty time slots at the end of each day are all accounted for when assigning the patients. 
Finally, constraints~\eqref{UB2:roS=1}--\eqref{UB2:roS=Ns} have the same meaning as~\eqref{UB2:roB=1}--\eqref{UB2:roB=Nb}, but operate on chairs.


\section{Computational results} 
\label{sec:results}

The solution approaches described in Section~\ref{sec:procedure} have been implemented in AMPL. The ILP models have been solved with Gurobi 9.5.1. Tests have been run on a Apple M1 MAX equipped with a 64 GB RAM, running under MacOS 13.4. 

We compare the results of our approaches, for both lower and upper bounds, with the ones obtained by Gurobi applied to the formulations. Different time limits have been set for the different steps of the approaches. 

In Section~\ref{subsec:case study} we describe the considered instances. 
Section~\ref{sec::results-P1} reports the comparison between formulations \eqref{P1} and $(AF_1)$. Section~\ref{sec::results-P2} reports the comparison between formulation~\eqref{SF2} and Procedure~\ref{proc:1}. 
The comparison between formulation~\eqref{SF3} and Procedure~\ref{proc:2} is reported in Section~\ref{sec::results-P3}. 
Finally, the upper bounds $(UB_1)$ and $(UB_2)$ are compared in Section~\ref{subsec::UB}, and the values of the three metrics are analyzed in Section~\ref{sec::overall}.


\subsection{Instances description} 
\label{subsec:case study}

To test our approach we use data provided by the outpatient cancer center of the San  Martino University Hospital in Genova (Italy). 

The center provides chemotherapy treatments for patients coming from the hematological and oncological departments of the hospital in an outpatient setting. 
In the center, there are 6 rooms available for the oncologist and hematologist visits, 26 chairs, and 27 beds for drug administration.   
The center is open from 8 to 17 from Monday to Friday, corresponding to $|H|= 54$ time slots, each 10 minutes long, available on each day of the planning horizon $T$, which is composed of five working days. 
However, patients can be visited only from 8 to 14, when the clinicians are available for visits and the rooms are open, corresponding to $|H_V|= 36$ time slots, 10 minutes long, for each day. 

Our approaches have been tested on a set of 52 instances representing the set of patients to be scheduled every week of the year 2019. 
In 2019 32,018 chemotherapy treatments have been recorded.
Patients are grouped into seven cancer pathologies, or cancer Macro-Groups (cMG), adapting the classification, reported by the American Joint Committee on Cancer (AJCC)~\citep{AJCC2017}. 
One group contains all the patients affected by hematologic cancers (HE) such as leukemia, lymphoma, and  myeloma, which represent about 30\% of the total number of patients treated. 
Among the oncological patients, the solid tumors have been classified into 6 categories mainly based on the specific body site or organ affected by cancer: Breast (BR), Gynecology (GY), Lung-respiratory (LU), Gastrointestinal (GI), Urology (UR), and Others (OT)~\citep{CLTT2022}. 
The visit time $v_p$ is assumed deterministic and it is 10 minutes, on average, for the oncological patients and 20 minutes for the hematologic patients, who usually require a more thorough examination by the clinician. 
The drug administration time $f_p$ ranges from one to four hours depending on the cMG.

In Table~\ref{tab:data1} the total number of patients $|P|$, the number and percentage of critical and non-critical patients, as well as  the number and percentage of patients for each cMG, and their average, minimum and maximum values are reported. 
The number of patients to be scheduled each week ranges from 430 to 719. 
The critical patients represent, on average, 28\% of the total number of accesses. 
With respect to the distribution of patients among the cMG, the most frequent are the hematology patients (HE) and breast cancer patients (BR) who account, together, for more than 50\% percent of the total number of treatments, in all the weeks, except two (week 33 and week 52). Indeed, in four out of 52 instances that percentage exceeds 70\%.
The data reported in Table~\ref{tab:data1} generally show a limited variation among weeks as for the distribution of patients in the cMGs, except for the hematologic case, with a Standard Deviation of about 25 patients, as appointments for chemotherapy treatments, in particular for solid tumors, are cyclically scheduled according to the treatment plan decided by the oncologist. 

\begin{table}[tbp]
  \centering
  \footnotesize
\tabcolsep 2.8pt
\begin{tabular}{l c r c r r r r r r r r r r r r r r}
\toprule
& 
& \multicolumn{2}{c}{\multirow{2}[-1]{*}{Critical}} & \multicolumn{14}{c}{Cancer Macro Groups} 
\\
\cmidrule(lr){5-18}          
&       & \multicolumn{2}{c}{patients} & \multicolumn{2}{c}{HE} & \multicolumn{2}{c}{GI} & \multicolumn{2}{c}{UR} & \multicolumn{2}{c}{GY} & \multicolumn{2}{c}{BR} & \multicolumn{2}{c}{OT} & \multicolumn{2}{c}{LU} 
\\
\cmidrule(lr){3-4}
\cmidrule(lr){5-6}
\cmidrule(lr){7-8}
\cmidrule(lr){9-10}
\cmidrule(lr){11-12}
\cmidrule(lr){13-14}
\cmidrule(lr){15-16}
\cmidrule(lr){17-18}
 & $|P|$ & \multicolumn{1}{c}{\#} & \multicolumn{1}{c}{\%} & \multicolumn{1}{c}{\#} & \multicolumn{1}{c}{\%} & \multicolumn{1}{c}{\#} & \multicolumn{1}{c}{\%} & \multicolumn{1}{c}{\#} & \multicolumn{1}{c}{\%} & \multicolumn{1}{c}{\#} & \multicolumn{1}{c}{\%} & \multicolumn{1}{c}{\#} & \multicolumn{1}{c}{\%} & \multicolumn{1}{c}{\#} & \multicolumn{1}{c}{\%} & \multicolumn{1}{c}{\#} & \multicolumn{1}{c}{\%} 
\\ 
\midrule
Total & 31942 & 7819  & 14.81 & 8872  & 16.80 & 3267  & 6.19 & 1909  & 3.62 & 1124  & 2.13 & 8459  & 16.02 & 4003  & 7.58 & 4308  & 8.16 
\\
    Avg & 614 & 150 & 28.48 & 171 & 32.31 & 63 & 11.90 & 37 & 6.95 & 22 & 4.09 & 163 & 30.81 & 77 & 14.58 & 83 & 15.69 \\
    Min   & 430   & 81    & 15.34 & 90    & 17.05 & 36    & 6.82 & 24    & 4.55 & 12    & 2.27 & 121   & 22.92 & 60    & 11.36 & 60    & 11.36 \\
    Max   & 719   & 196   & 37.12 & 214   & 40.53 & 81    & 15.34 & 51    & 9.66 & 37    & 7.01 & 194   & 36.74 & 94    & 17.80 & 115   & 21.78 \\
    \bottomrule
    \end{tabular}%
\caption{Instances patient data.}
\label{tab:data1}%
\end{table}%

The MCP gives the assignment of cMGs to the days of the week and exam rooms available in the outpatient cancer center. In the considered cancer center, the weekly pattern is defined at the beginning of each month and it is cyclically repeated over the month. 
The MCP is designed as described in~\cite{CLTT2022}, based on the demands for chemotherapy treatments for each cMG which come from the clinical treatment plans and on the clinicians' availability and skills. According to the cancer center policy, three exam rooms are assigned every working day to the hematologic patients, the other three rooms are instead assigned to oncological patients.


\subsection{Results on problem $(P_1)$}
\label{sec::results-P1}

We compare the performance of the formulation \eqref{P1}  and  the approach based on the aggregate formulation, that first solves $(AF_1)$ and then applies Algorithms~\ref{alg:1}, \ref{alg:2} and \ref{alg:3}.
The time limit has been set to 300 seconds for each run of the solver. 
Aggregated results are shown in Table~\ref{tab:P1-comparison}, which reports the number of instances solved to optimality, the number of feasible solutions found, optimality gaps, and CPU times. 
The aggregate formulation $(AF_1)$ significantly outperforms~\eqref{P1}. 
In fact, it solves all the instances in less than one minute, while \eqref{P1} is able to solve to optimality only two instances and cannot find a feasible solution in five. 
Therefore, we will use the aggregate formulations also for solving problems $(P_2)$ and $(P_3)$.

\begin{table}[htbp]
\centering
\begin{tabular}{l c c} 
\toprule
  & $(F_1)$   & $(AF_1)$+Algorithms 1, 2, 3 \\
\hline
\# time limit without optimal solution (out of 52) & 50    & 0 \\
\# optimal solution found (out of 52) & 2     & 52 \\
    \# LB found (out of 52) & 47    & 52 \\
    Average percentage optimality gap [\%]  & 19.14 & 0.00 \\
    Average CPU time for all instances [s] & TL & 17.86 \\ 
    Maximum CPU time for all instances [s] & TL & 57.88 \\
    \midrule
    Average $(AF_1)$ CPU time  [s] &       & 11.20 \\
    Maximum $(AF_1)$ CPU time [s] &       & 49.95 \\
    Average algorithms CPU time [s] &       & 6.66 \\
    Maximum algorithms CPU time [s] &       & 8.49 \\
\bottomrule
    \end{tabular}%
      \caption{Problem $(P_1)$: comparison between $(F_1)$ and the proposed approach ($(AF_1)$+Algorithms 1, 2, 3).}
  \label{tab:P1-comparison}%
\end{table}%


\subsection{Results on problem $(P_2)$}
\label{sec::results-P2}

Although the aggregate formulation outperforms \eqref{P1} for problem $(P_1)$, it is time consuming for problem $(P_2)$. 
Thus, we developed Procedure~\ref{proc:1} (see Section~\ref{s:solvingP2}): we first restrict the problem according to the solution of problem $(P_1)$ and then decompose it into a set of daily problems. 
The resulting solution is then used as a warm start for $(AF_2)$. 
We now compare the two approaches: the $(AF_2)$ formulation solved by the solver as is, and Procedure~\ref{proc:1}. 
A 600 seconds time limit is set for the solution of $(AF_2)$, a 60 seconds time limit is set for the solution of every single day problem, and 300 seconds are assigned to $(AF_2)$ starting from the warm start.

Results are summarized in Table~\ref{tab:P2-comparison}: for both approaches, we report the number of instances 
in which the optimal solution is found, those where an upper bound is found, the average and maximum CPU time when an optimal solution is found, and the average CPU time for all the instances. 
Results clearly show that Procedure~\ref{proc:1} outperforms $(AF_2)$. 
In fact, $(AF_2)$ is able to find the optimal solution only in 9 instances, while Procedure~\ref{proc:1} solves to optimality 46 instances out of 52. 
When $(AF_2)$ is not able to find the optimal solution, Procedure~\ref{proc:1} improves significantly upon $(AF_2)$, as it is able to obtain much shorter waiting times: it decreases by about 12 time slots (about 2 hours) on average, and it may reduce up to 154. 
Further, $(AF_2)$ fails in finding a feasible solution in five instances, while Procedure~\ref{proc:1} finds a feasible solution in all the instances.
In the 6 instances for which optimality is not proved, even the single-day problems cannot be solved to optimality. 
Further, the upper bound provided by Procedure~\ref{proc:1} is very close to zero: it is equal to one in four instances, and to two and three in the remaining two, respectively. 
The solution provided by the first step (restricted and decomposed problem) is usually very good and the second step is able to improve it only in one instance. 
As a consequence, solving $(AF_2)$ starting from the warm start is usually very fast. 
Procedure~\ref{proc:1} is about 36 times faster than $(AF_2)$ on average. 

\begin{table}[htbp]
\centering
\begin{tabular}{@{}l c c c c@{}}
\toprule
&  & \multicolumn{3}{c}{Procedure~\ref{proc:1}} 
\\
\cmidrule(lr){3-5}
    &       & {Restricted and} & $(AF_2)$ with & {overall} 
\\
     & $(AF_2)$  & decomposed  $(P_2)$  & warm start &  approach 
\\
\midrule
  \# optimal solution found (out of 52) & 9      & 46 & 46 & 46\\
 \# UB found (out of 52) & 47    &    52 & 52 & 52 
\\
Avg CPU time with optimal solution  [s] & 170.73 & 17.64 & 0.90 & 18.54
\\
Max CPU time with optimal solution [s] & 291.85 & 167.49 & 1.23  & 168.72
\\ 
Avg CPU time for all the instances  [s] & 535.41 & 32.17 & 35.38 & 67.55\\
\bottomrule
\end{tabular}%
\caption{Problem $(P_2)$: comparison between $(AF_2)$ and Procedure~\ref{proc:1}.}
\label{tab:P2-comparison}%
\end{table}%


\subsection{Results on problem $(P_3)$}
\label{sec::results-P3} 

We now compare the results on problem ($P_3$). 
In Table~\ref{tab:P3-LB-comparison} the results obtained solving \eqref{SF3} within a 600 seconds time limit are compared with the results of Procedure~\ref{proc:2}. 
For the latter, 90 seconds are assigned to each daily problem, 60 seconds to each iteration of the $k$-opt local search, and a 600 seconds time limit is set for the overall heuristic procedure. 
We consider all the instances: for those for which the optimal solution of ($P_2$) is not known, we consider the best solution found to set the value of $\bar{v}^t_2$. 
Using the solution of Procedure~\ref{proc:2} as a warm start does not produce improvement, therefore we do not apply the step. We use, instead, the bounds described in Section~\ref{sec::upperBounds} as a comparison. Table~\ref{tab:P3-LB-comparison} shows that Procedure~\ref{proc:2} clearly outperforms \eqref{SF3} in terms of the number of best solutions found and the gap with respect to the best lower and upper bound. 
In fact, Procedure~\ref{proc:2} finds the best solution in 44 instances, while \eqref{SF3} is able to find only 10 best lower bounds. 
The gap w.r.t. the best lower bound known is 0.11\% on average for Procedure~\ref{proc:2} and it may rise up to 1.67\%; the optimality gap is 6.86\% on average and may rise up to 10.82\%. 
The gaps provided by \eqref{SF3} are far worse: about 2\% and 9\% for the best lower bound gap, and about 8.8\% and 18.7\% for the optimality gap. Further, \eqref{SF3} is far slower than the heuristic approach. 
Although even the first heuristic step outperforms  the formulation, the $k$-opt step is able to further improve the quality of the obtained solutions.

\begin{table}[htbp]
\centering
\begin{tabular}{@{}lcccc@{}}
\toprule
    &  & \multicolumn{3}{c}{Procedure~\ref{proc:2}} 
\\ 
\cmidrule(lr){3-5}
         &       & {Restricted and} & \multirow{2}{*}{$k$-opt} & {overall } \\
        & $(AF_3)$      & {decomposed $(P_3)$} & {} & {approach} \\
\midrule
Average CPU time [s] & TL    & 64.61 & 149.46 & 248.06 \\

Maximum CPU time [s] & TL    & 141.56 & 361.75 & 435.85 \\

    \# best LB found (out of 52) & 10    & 16    & 44    & - \\
    
Average percentage gap w.r.t. best LB [\%] & 2.19 & 1.05 & 0.11 & - \\
    
Maximum percentage gap w.r.t. best LB [\%] & 9.15 & 5.50 & 1.67 & - \\

Average percentage gap w.r.t. best UB [\%] & 8.80 & 7.74 & 6.86 & - \\
    
Maximum percentage gap w.r.t. best UB [\%] & 18.70 & 11.34 & 10.82 & - \\
\bottomrule
\end{tabular}%
\caption{Problem $(P_3)$: comparison between $(AF_3)$ and Procedure~\ref{proc:2}.}
\label{tab:P3-LB-comparison}%
\end{table}%

\begin{figure}[htbp]
\centering
\includegraphics[width=0.65\textwidth]{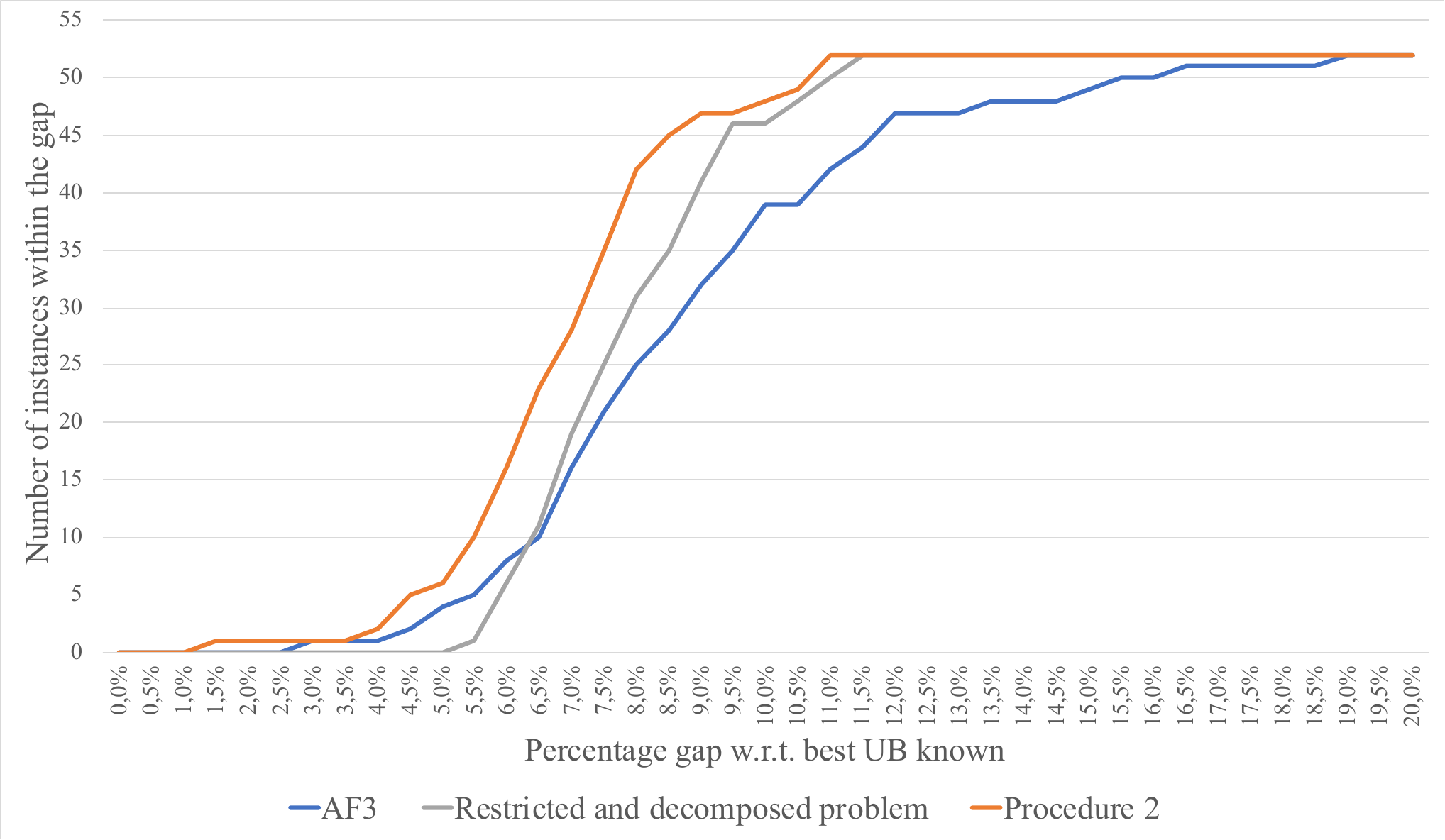}
\caption{$(P_3)$ performance profile: gap w.r.t. the best upper bound known.}
\label{fig:pp-bestUBound}
\end{figure}

The quality of the heuristic solutions is further highlighted in Figure~\ref{fig:pp-bestUBound}, 
where the performance profiles related to the gap w.r.t. upper bound known are represented for~\eqref{SF3} and for the restricted and decomposed problem and overall Procedure~\ref{proc:2}. The figure shows that our approach, even without the $k$-opt neighborhood search, outperforms the gap provided by \eqref{SF3}. 

\subsection{Quality of the upper bound for ($P_3$)}
\label{subsec::UB} 

To assess the quality of the upper bound for ($P_3$) (see Section~\ref{sec::upperBounds}) we compare the upper bound computed by \eqref{SF3}, a straightforward upper bound (the number of non-critical patients), the multi knapsack based upper bound ($UB_1$), and the improved upper bound ($UB_2$). A 600 seconds time limit is assigned to \eqref{SF3} and a 60 seconds time limit to $(UB_2)$. 
Table~\ref{tab:P3-UB-comparison} shows that the best results are obtained by $(UB_2)$: it finds the best upper bound in 48 instances out of 52 and has a very small gap w.r.t. the best upper bound in the remaining four - on average 0.8\% and never above 2.5\%.  \eqref{SF3} is the second best approach in terms of quality of the obtained upper bound: however, its gap w.r.t. the best upper bound may rise up to about 7\% and it reaches the time limit in all the instances, thus requiring about 10 times the CPU time required by $(UB_2)$. 
In fact, $(UB_2)$ reaches the time limit only in 20 instances out of 52 and is given only 60 seconds as a time limit. The other two upper bounds, despite being very fast, provide far worse results. In fact, their bounds have a gap w.r.t. the best upper bound known of about 6\% on average and up to more than 20\% in the worst case. As a consequence, \eqref{SF3} cannot properly assess the quality of its own feasible solutions if it relies on its own lower bound. Indeed, the quality of ($UB_2$) makes it possible to significantly reduce the optimality gap.

\begin{table}[htbp]
\centering
\begin{tabular}{lcccc} 
\toprule
                & \multirow{2}{*}{$(AF_3)$} & \# of non-critical  & \multirow{2}{*}{$(UB_1)$}    & \multirow{2}{*}{$(UB_2)$} \\ 
                &  &  patients &     &  \\
\midrule
    \# best UB found &        12    & 9     & 9     & 48 \\
Average percentage gap w.r.t. best UB [\%]      & 3.61 & 6.00 & 5.69 & 0.08 
\\
Maximum percentage gap w.r.t. best UB [\%]      & 6.97 & 31.52 & 24.90 & 2.47 
\\
Average CPU time [s]        & TL & -- & 0.22 & 33.99 
\\
Maximum CPU time [s]       & TL & -- & 0.48 & TL \\
\bottomrule
\end{tabular}%
\caption{Comparison of upper bounds.}
\label{tab:P3-UB-comparison}%
\end{table}


\subsection{Metrics}
\label{sec::overall}

Let us now focus on the three considered metrics: the number of scheduled patients~\eqref{obf:max_p_surr}, the waiting time~\eqref{OBJ:waitingTime-xyz}, and the number of non-critical patients who receive the infusion on a chair~\eqref{OBJ:chairPreferences-xyz}. 

As mentioned, the instances represent 52 weeks of a year. 
It may happen, therefore, that in some weeks not all the days are working days, due to holidays and vacations: it is the case for 6 of the considered weeks. 
In Table~\ref{tab:treated-patients}, only the 46 weeks in which there are five working days are considered. The table reports the number of weeks in which the patients cannot be all scheduled, and their average and maximum percentage, for the overall set of patients, for each pathology, and for the critical and non-critical patients. 
Although in 25 weeks there are some unscheduled patients, their number is low, less than 1.5\% on average (considering only the 25 weeks where not all the patients can be scheduled) and never above 3.5\%. However, despite being in general acceptable, the number of unscheduled patients changes significantly from one pathology to another. For some pathologies, the average percentage of unscheduled patients is high, for instance, it is about 18\% in a single week for gynecology patients, and about 14\% over 16 weeks for urology patients. 
Hematologic patients, instead, are always all scheduled in the requested week, despite being the largest group of patients (about 32\% of the patients), as the largest number of visit room slots in the MCP is assigned to them. 
Breast cancer patients are the second largest group (about 30\% of the total number of requests) but are usually assigned about one-third of the visit room slots w.r.t. those assigned to the hematology patients. 
Even if hematology patients require twice the visit time as the others, it is nevertheless an unbalanced resource assignment, resulting in 6 weeks where on average about 3.6\% of breast cancer patients are not scheduled. 
Further, when less than five working days are available, some groups of patients cannot be scheduled at all, as they are assigned visit rooms on only one day of the week, based on the MCP, and on this day the center is not open. 
Critical patients seem to be more difficult to handle then non-critical ones. 
In fact, despite being on average less than 30\%, they cannot be all scheduled in 25 weeks, while the non-critical ones cannot be all scheduled only in 17. Further, the average percentage of untreated patients among the critical ones is more than twice the percentage of unscheduled patients among the non-critical ones.

In the optimal solutions of problem ($P_1$) the values of the other two metrics are usually poor. 
As for the waiting time metrics, it is strictly positive in all the weeks, it is on average about 30 hours, and may rise up to 35. Instead, it is dramatically reduced by solving ($P_2$). 
In fact, as discussed in Section~\ref{sec::results-P2}, in the solution of ($P_2$) the waiting time is equal to 0 in all the weeks except for 6 weeks, where it is equal to half an hour, at most. 
We can therefore infer that the two metrics are not related and that there are many optimal solutions of ($P_1$) among which we can select the ones with the minimum waiting time. 
A similar remark applies also to the third metrics, i.e., the number of non-critical patients receiving the infusion on a chair. 
In fact, solving ($P_3$) produces an improvement in the third metrics of about 45\% with respect to both the solutions of ($P_1$) and ($P_2$), selecting among equivalent optima the one providing significantly better values of the third metrics.

\begin{table}[htbp]
\small
\centering
\tabcolsep 3pt
\begin{tabular}{@{}lc|ccccccc|cc@{}}
\toprule
     & & \multicolumn{7}{|c|}{Cancer Macro Groups}  &   &
\\
\cmidrule(lr){3-9}
& overall & HE    & GI    & UR    & GY    & BR    & OT    & LU     & critical & non-critical 
\\
\midrule
\# weeks with & & & & & & & & & 
\\
unscheduled patients & 25 & 0    & 6    & 16    & 1    & 6    & 2    & 4    & 25 & 17
\\
\midrule
Average percentage  & & & & & & & & &
\\
 unscheduled patients [\%] & 1.36 & 0.00 & 5.99 & 14.05 & 18.92 & 3.61 & 9.34 & 7.24 & 2.80 & 1.23
\\
\midrule
Maximum percentage & & & & & & & & & 
\\
unscheduled patients [\%] & 3.33 & 0.00 & 11.11 & 29.41 & 18.92 & 7.22 & 12.20 & 15.29 & 8.06 & 2.18
\\
\bottomrule
\end{tabular}%
\caption{Unscheduled appointments analysis.}
\label{tab:treated-patients}%
\end{table}%


\section{Conclusions}
\label{sec:conclusions}

We addressed a chemotherapy treatment planning and multi-appointment scheduling problem in a shared cancer center. We consider three different metrics associated with the patient's comfort. 
We tackle the resulting multi-objective problem by solving three problems in sequence, as the three objectives are in a hierarchical order. 
Alternative formulations, heuristic approaches, and bounds have been developed for the three problems.  The resulting approaches outperform the state-of-the-art solvers and allow us to obtain very good quality solutions in reasonable computational times.

As for future development, we remark that, as far as the first metrics is concerned, the allocated resources cannot always allow to treat all the patients. Although, in the current practice, the issue can be dealt with by allowing a little amount of overtime, it is worth exploiting the proposed approach to evaluate the suitability of the allocated resources in other scenarios and investigate alternative policies and resource dimensioning. 
As future research direction, we will address the problem of jointly designing the MCP, planning the appointments, and scheduling the patients' activities. 
Further, we will address the case where the bad health condition of the patients may prevent to receive the treatment and embed the related uncertainty in the problem.


\bibliographystyle{apalike}

\bibliography{references2}

\end{document}